\documentclass[10pt]{amsart}
\usepackage{amssymb}

\theoremstyle{definition}
\newtheorem{thm}{Theorem}[section]
\newtheorem{lem}[thm]{Lemma}
\newtheorem{prp}[thm]{Proposition}
\newtheorem{dfn}[thm]{Definition}
\newtheorem{cor}[thm]{Corollary}

\newtheorem{exa}[thm]{Example}
\newtheorem{wrn}[thm]{Warning}
\newenvironment{pff}{{\emph{Proof:}}}{\QED}

\newcommand{\bit}{\begin{itemize}}
\newcommand{\eit}{\end{itemize}}
\newcommand{\beq}{\begin{equation}}
\newcommand{\eeq}{\end{equation}}
\newcommand{\beqr}{\begin{eqnarray*}}
\newcommand{\eeqr}{\end{eqnarray*}}
\newcommand{\bal}{\begin{align*}}
\newcommand{\eal}{\end{align*}}
\newcommand{\ts}{\textstyle}
\newcommand{\rsz}[1]{\raisebox{0ex}[0.8ex][0.8ex]{$#1$}}

\newcommand{\af}{\alpha}
\newcommand{\bt}{\beta}

\newcommand{\ep}{\varepsilon}

\newcommand{\et}{\eta}

\newcommand{\sm}{\sigma}

\newcommand{\ph}{\varphi}
\newcommand{\ps}{\psi}
\newcommand{\rh}{\rho}

\newcommand{\ta}{\tau}

\newcommand{\Gm}{\Gamma}

\newcommand{\Z}{{\mathbf{Z}}}
\newcommand{\R}{{\mathbf{R}}}
\newcommand{\C}{{\mathbf{C}}}
\newcommand{\N}{{\mathbf{N}}}

\newcommand{\Ker}{{\mathrm{Ker}}}
\newcommand{\inv}{{\mathrm{inv}}}
\newcommand{\tsr}{{\mathrm{tsr}}}

\newcommand{\id}{{\mathrm{id}}}
\newcommand{\ev}{{\mathrm{ev}}}
\newcommand{\sint}{{\mathrm{int}}}

\newcommand{\dist}{{\mathrm{dist}}}

\newcommand{\rank}{{\mathrm{rank}}}
\newcommand{\GL}{{\mathrm{GL}}}

\newcommand{\op}{{\mathrm{op}}}
\newcommand{\Prim}{{\mathrm{Prim}}}

\newcommand{\dirlim}{\displaystyle \lim_{\longrightarrow}}
\newcommand{\invlim}{\displaystyle \lim_{\longleftarrow}}
\newcommand{\Mi}{M_{\infty}}

\newcommand{\andeqn}{\,\,\,\,\,\, {\text{and}} \,\,\,\,\,\,}
\newcommand{\QED}{\rule{0.4em}{2ex}}

\newcommand{\ca}{C*-algebra}
\newcommand{\ct}{continuous}
\newcommand{\pj}{projection}

\newcommand{\nbhd}{neighborhood}
\newcommand{\cpt}{compact Hausdorff}
\newcommand{\hm}{homomorphism}

\newcommand{\Wolog}{Without loss of generality}
\newcommand{\Tfae}{The following are equivalent}
\newcommand{\tfae}{the following are equivalent}
\newcommand{\ifo}{if and only if}
\newcommand{\rsha}{recursive subhomogeneous algebra}
\newcommand{\Rsha}{Recursive subhomogeneous algebra}

\newcommand{\rshd}{recursive subhomogeneous decomposition}
\newcommand{\mops}{mutually orthogonal \pj s}
\newcommand{\tdim}{topological dimension}
\newcommand{\mms}{maximum matrix size}

\title[Recursive subhomogeneous algebras]{Recursive subhomogeneous
algebras}

\author{N.\  Christopher Phillips}

\address{Department of Mathematics, University  of Oregon,
       Eugene OR 97403-1222, USA}

\thanks{Research
      partially supported by NSF grants DMS 9400904 and DMS 9706850.}

\begin{document}

\setcounter{section}{-1}

\begin{abstract}
We introduce and characterize a particularly tractable class of unital
type~1 \ca s with bounded dimension of irreducible representations.
Algebras in this class are called \rsha s, and they have an inductive
description (through iterated pullbacks) which allows one to carry over
from algebras of the form $C (X, M_n)$ many of the constructions
relevant in the study of the stable rank and K-theory of simple
direct limits of homogeneous \ca s.
Our characterization implies in particular that
if $A$ is a separable \ca\  whose irreducible representations all
have dimension at most $N < \infty$, and if for each $n$ the
space of $n$-dimensional irreducible representations has finite
covering dimension, then $A$ is a \rsha.
We demonstrate the good properties of this class by proving
subprojection and cancellation theorems in it.

Consequences for simple direct limits of \rsha s, with applications
to the transformation group \ca s of minimal homeomorphisms,
will be given in a separate paper.
\end{abstract}

\maketitle

\section{Introduction}

In recent years, a number of results have been proved about the
stable and real rank and the unstable K-theory of certain kinds of
direct limits of homogeneous \ca s, usually assuming some kind
of slow dimension growth and sometimes assuming simplicity.
See for example \cite{DNNP}, \cite{BDR}, \cite{Bl4},
\cite{MP}, and \cite{Gd}.
The subprojection and cancellation theorems for algebras of the
form $C (X, M_n)$ are important ingredients in many of the proofs.
These theorems essentially say that if $\rank (p) - \rank (q)$
is large enough compared to the dimension of $X$, then $q$ is a
sub\pj\  of $p$, and that if $\rank (p)$
is large enough compared to the dimension of $X$, and if
$p \oplus e$ is Murray-von Neumann equivalent to $q \oplus e$,
then $p$ is Murray-von Neumann equivalent to $q$.
See Chapter~9 of \cite{Hs} for the original formulation in terms of
vector bundles over finite complexes, and see Theorem~2.5 of \cite{Gd},
Lemma~3.4 of \cite{MP}, and Lemma~1.5 of \cite{Ph1} for
\ca ic formulations.

In this paper, we introduce a class of type~1 \ca s, much more general
than those of the form $C (X, M_n)$, in which analogs of these theorems
are still true.
Algebras in this class, called \rsha s, can be expressed as iterated
pullbacks of algebras of the form $C (X, M_n)$, in such a way that an
inductive argument reduces these theorems to ``relative'' versions
of the same theorems for $C (X, M_n)$.
In a second paper \cite{PhX},
we use these results to generalize some of the known
results on direct limits of direct sums of homogeneous \ca s to
direct limits of \rsha s.
Combining these results on direct limits with the work of
Qing Lin \cite{Ln}, we will obtain information on the order on
the $K_0$ groups of \ca s of minimal homeomorphisms.

\Rsha s include a number of algebras which have already played
significant roles in the study of direct limits of type~1 \ca s,
and also elsewhere.
The following are all \rsha s:
\begin{itemize}
\vspace{-0.5ex}
\item
Finite direct sums of algebras of the form $C (X, M_n)$ for
$X$ \cpt\  and $n \geq 1$.
\item
Dimension drop intervals (as used in \cite{El})
and matrix algebras over them.
\item
The building blocks used in the classification theorems of
\cite{JS}, \cite{Th}, and \cite{Mg}.
\item
The noncommutative CW complexes of \cite{Pd}.
\item
Section algebras of locally trivial continuous
fields over $X$ with fiber $M_n$ (with arbitrary Dixmier-Douady class).
\item
The algebras $A_Y$ arising in Qing Lin's study \cite{Ln} of the
transformation group \ca s of minimal homeomorphisms, provided
$\sint (Y) \neq \varnothing$.
(See Example~\ref{A6} for more details.)
\end{itemize}

It is obvious from the definition (given in Section~1) that \rsha s
are unital type~1 \ca s with a finite upper bound on the possible
dimensions of irreducible representations.
Among such algebras, they at first seem rather special, but in fact
there are a great many \rsha s.
We prove in Section~2 a characterization theorem which implies the
following:
If $A$ is a separable \ca\  whose irreducible representations all
have dimension at most $N < \infty$, and if for each $n$ the
space of $n$-dimensional irreducible representations has finite
covering dimension, then $A$ is a \rsha.

Work related to our first two sections has been done in \cite{Vs}.
(We are grateful to George Elliott for pointing out this reference.)
Specifically, Theorem~4 of \cite{Vs} gives a decomposition of an
arbitrary subhomogeneous \ca\  which is related to that in our
Proposition~\ref{N11}.
(Note that we show the decomposition described in Proposition~\ref{N11}
holds under the hypotheses of Theorems~\ref{N13} and~\ref{N14}.)
However, the description of \cite{Vs} is much more complicated and
difficult to work with.
In our applications (see, for example, \cite{PhX}),
the spaces of irreducible representations of each
fixed dimension will all have finite covering dimension, in which case
Theorem~\ref{N14} guarantees the applicability of Proposition~\ref{N11}.
Moreover, in the work outlined in Section~6 of \cite{LP},
we apparently need \rshd s satisfying an additional condition
(a suitable a priori bound on the ``strong covering number'').
As far as we know, the decomposition of \cite{Vs}
need not satisfy any analog of this condition;
nor do we know how to obtain it
for the decomposition of Proposition~\ref{N11}.

There are four sections.
The first section contains the definition of a \rsha, examples,
elementary closure properties of the class, and related results.
The second gives a characterization of separable \rsha s, both in
general and under the condition that the spaces of
$n$-dimensional irreducible representations have finite
covering dimension.
In the third section, these results are used to give stronger
closure properties of the class of \rsha s.
In particular, quotients, unitized ideals, suitable pullbacks, and
tensor products of \rsha s are all again \rsha s.
However, it is shown by example that subalgebras of \rsha s need not
be \rsha s.
In Section~4 we prove the sub\pj\  and cancellation theorems for
\rsha s, and an analogous result about unitaries.
Applying the results of the previous sections, we obtain
subprojection and cancellation theorems for separable
\ca s with a finite upper
bound on the possible dimensions of irreducible representations,
under suitable assumptions
on the dimensions of subspaces of the primitive ideal space.
We also give a sufficient condition for the $K_1$-group to be
determined by the unitaries in the algebra, without using matrices.

In this paper,
we let $\dim (X)$ denote the {\emph{modified}} covering dimension
of the topological space $X$.
See Definition~10.1.3 of \cite{Pr}.
This is the same as the usual covering dimension
(Definition~3.1.1 of \cite{Pr}) on compact Hausdorff spaces.
Since we only use other spaces in Section~2, we postpone further
discussion to there.
(Warning on terminology: ``bicompact'' in  \cite{Pr} is what is usually
called compact Hausdorff, and ``Tihonov'' in  \cite{Pr}
is what is usually called
completely regular, that is, points are closed and can be separated
from closed sets by continuous functions.)

I am grateful to Larry Brown, Maurice Dupr\'{e}, Qing Lin,
Gert Pedersen, and Claude Schochet
for useful conversations and email correspondence.
In particular, the definition of a \rsha\  arose from an effort with
Qing Lin to impose structure on the algebras arising in his
study \cite{Ln} of the
transformation group \ca s of minimal homeomorphisms
(see Example~\ref{A6}), and
Larry Brown suggested the possible relevance of the finite type
condition used in Section~2.
This work was carried out during a sabbatical year at
Purdue University, and I am grateful to that institution
for its hospitality.

Some of the results of this paper were announced in \cite{LP}.

\section{\Rsha s}

In this section, we give the definition of a \rsha, along with several
remarks and some useful terminology.
We then give a number of examples, including all those mentioned in
the introduction.
After that, we give some closure properties of the class of \rsha s
which can be proved easily from the definition.
In particular, direct sums and corners of \rsha s are \rsha s, and
tensor products of \rsha s with $C (X)$ and $M_n$ are \rsha s.
Some stronger closure properties will be obtained in Section~3.

\begin{dfn}\label{A1}
A {\emph{\rsha}} is a \ca\  given by the following recursive
definition.
\bit
\item[(1)]
If $X$ is a \cpt\  space and $n \geq 1$, then
$C(X, M_n)$ is a \rsha.
\item[(2)]
If $A$ is a \rsha, $X$ is a \cpt\  space,
$X^{(0)} \subset X$ is closed,
$\ph \colon A \to C \left( \rsz{ X^{(0)} }, \, M_n \right)$ is any
unital \hm, and
$\rh \colon C (X, M_n) \to C \left( \rsz{ X^{(0)} }, \, M_n \right)$
is the restriction \hm, then the pullback
\[
A \oplus_{C ( X^{(0)}, \, M_n )} C (X, M_n) =
  \{ (a, f) \in A \oplus C (X, M_n) \colon \ph (a) = \rh (f) \}
\]
(compare with Definition~2.1 of \cite{Pd}) is a \rsha.
\eit

Note that in~(2) the choice $X^{(0)} = \varnothing$ is allowed
(in which case $\ph = 0$ is allowed).
Thus the pullback could be an ordinary direct sum.

It is convenient in several situations (such as consideration of
corners) to allow the zero algebra to be a \rsha.
\end{dfn}

Examples will be presented shortly, but some terminology and a
warning are appropriate first.

\begin{dfn}\label{A2}
We adopt the following standard notation for \rsha s.
{}From the definition, it is clear that any \rsha\  can be written in
the form
\[
R = \left[ \cdots \rule{0em}{3ex} \left[ \left[
  C_0 \oplus_{C_1^{(0)}} C_1 \right]
 \oplus_{C_2^{(0)}} C_2 \right] \cdots \right]
            \oplus_{C_l^{(0)}} C_l,
\]
with $C_k = C (X_k, \, M_{n (k)})$ for \cpt\  spaces $X_k$ and positive
integers $n (k)$, with
$C_k^{(0)} = C {\ts{ \left( \rsz{ X_k^{(0)} }, \, M_{n (k)} \right) }}$
for compact
subsets $X_k^{(0)} \subset X_k$ (possibly empty), and where the maps
$C_k \to C_k^{(0)}$ are always the restriction maps.
An expression of this type will be referred to as a
{\emph{decomposition}} of $R$, and the notation used here
will be referred to as the standard notation for a decomposition.

Associated with this decomposition are:
\bit
\item[(1)]
its {\emph{length}} $l;$
\item[(2)]
the {\emph{$k$-th stage algebra}}
\[
R^{(k)} = \left[ \cdots \rule{0em}{3ex} \left[ \left[
  C_0 \oplus_{C_1^{(0)}} C_1 \right]
 \oplus_{C_2^{(0)}} C_2 \right] \cdots \right]
            \oplus_{C_k^{(0)}} C_k,
\]
obtained by using only the first $k + 1$ algebras
$C_0, C_1, \dots, C_k$;
\item[(3)]
its {\emph{base spaces}} $X_0, X_1, \dots, X_l$ and
{\emph{total space}} $X = \coprod_{k = 0}^l X_k;$
\item[(4)]
its {\emph{matrix sizes}} $n (0), \dots, n (l)$, and
{\emph{matrix size function}} $m \colon X \to \N \cup \{0\}$,
defined by $m (x) = n (k)$ when $x \in X_k$ (this is called
the {\emph{matrix size of $A$ at $x$}});
\item[(5)]
its {\emph{minimum matrix size}} $\min_k n (k)$ and
{\emph{maximum matrix size}} $\max_k n (k);$
\item[(6)]
its {\emph{topological dimension}} $\dim (X)$
(the covering dimension of $X$ \cite{Pr}, Definition~3.1.1;
here equal to  $\max_k \dim (X_k)$), and
{\emph{topological dimension function}} $d \colon X \to \N \cup \{0\}$,
defined by $d (x) = \dim (X_k)$ when $x \in X_k$ (this is called
the {\emph{topological dimension of $A$ at $x$}});
\item[(7)]
its {\emph{standard representation}}
$\sm = \sm_R \colon R \to \bigoplus_{k = 0}^l C (X_k, \, M_{n (k)} )$,
defined by forgetting the restriction to a subalgebra in each of
the fibered products in the decomposition;
\item[(8)]
the associated
{\emph{evaluation maps}} $\ev_x \colon R \to M_{n (k)}$ for $x \in X_k$,
defined to be the restriction of the usual evaluation map to
$R$, identified with a subalgebra of
$\bigoplus_{k = 0}^l C (X_k, \, M_{n (k)} )$ via $\sm$.
\eit
\end{dfn}

\begin{wrn}\label{A3}
The decomposition of a \rsha\  is highly nonunique (as will become
clear from some of the examples).
Throughout this paper, we will tacitly assume (unless otherwise
specified) that every \rsha\  comes given with some decomposition.
In particular, we will refer to the length, matrix sizes, etc.\ %
of a \rsha, by which we mean the corresponding quantities for
a tacitly understood given decomposition of the algebra.
\end{wrn}

We now give examples; for later use, some are given as parts
of propositions.

\begin{exa}\label{A4}
Any finite direct sum of \ca s of the form $C (X, M_n)$ (trivial
homogeneous algebras) is a \rsha.
More generally, any finite direct sum of \rsha s is a \rsha.
\end{exa}

\begin{exa}\label{A4.3}
The noncommutative CW-complexes of Section~11 of \cite{Pd}
are \rsha s (with particularly nice base spaces) whenever they are
unital.
In particular, the ``dimension drop intervals''
\[
\widetilde{I}_n =
   \{ f \in C ([0, 1], \, M_n) \colon f (0), \, f (1) \in \C \cdot 1 \},
\]
and the matrix algebras $M_k \left( \rsz{ \widetilde{I}_n } \right)$,
are \rsha s.
(These were introduced in \cite{BEEK},
and have since played a
significant role in the stably finite classification theory;
see, for example, \cite{El} and \cite{DG}.)

To write $\widetilde{I}_n$ as a \rsha\  of length $1$,
set
\[
C_0 = \C \oplus \C, \,\,\,\,\,\, C_1 = C ([0, 1], \, M_n), \andeqn
C_1^{(0)} = M_n \oplus M_n,
\]
and define maps from $C_0$ and $C_1$ to $C_1^{(0)}$ by
$\ph_1 (\af, \bt) = (\af \cdot 1, \, \bt \cdot 1)$ and
$\rh_1 (f) = (f (0), \, f (1))$.
To get $M_k \left( \rsz{ \widetilde{I}_n } \right)$,
form $k \times k$ matrices over everything.
(Compare with Lemma~\ref{A7} below.)

Algebras which are generalizations of the dimension drop intervals,
but are still one dimensional noncommutative CW-complexes, have been
used as building blocks in classification theorems for simple direct
limits in \cite{JS}, \cite{Th}, and \cite{Mg}.
\end{exa}

\begin{exa}\label{A6}
Let $X$ be an infinite compact metric space,
let $h$ be a minimal homeomorphism
of $X$, and let $Y \subset X$ be closed.
Let $A = C^* (\Z, X, h)$ be the transformation group \ca,
with unitary $u$ representing the generator of $\Z$.
Let
\[
A_Y = C^* ( C (X), \, u C_0 (X \setminus Y)) \subset C^* (\Z, X, h).
\]
The computations in the proof of Theorem~3 of \cite{Ln}
(also see Section~3 of \cite{LP}) show that,
if $\sint (Y) \neq \varnothing$, then $A_Y$ is a \rsha, with
topological dimension at most $\dim (X)$.
(This uses the fact that if $Z \subset X$ is closed, then
$\dim (Z) \leq \dim (X)$. See Proposition~3.1.5 of \cite{Pr}.)

It follows immediately that if $y_0 \in X$, then
$A_{\{y_0\}}$ is a direct limit of \rsha s, with all maps in
the direct system injective and unital.
In fact, if
\[
Y_1 \supset Y_2 \supset \cdots, \,\,\,\,\,\, \sint (Y_n) \neq \varnothing,
   \andeqn \bigcap_{n = 1}^{\infty} Y_n = \{y_0\},
\]
then $A_{\{y_0\}} = \overline{ \bigcup_{n = 1}^{\infty} A_{Y_n} }$.
\end{exa}

We next show that locally trivial continuous fields are \rsha s.
This follows from the following more general result, which
has other useful consequences.

\begin{prp}\label{A5}
Let $X$ be a \cpt\  space, let $E$ be a locally trivial continuous
field over $X$ with fiber $M_n$, and let $B = \Gm (E)$ be the
corresponding section algebra.
(We make no restriction on the Dixmier-Douady class.)
Let $A$ be a \ca\  with a \rshd\  of length $l_0$, \tdim\  $d_0$,
and \mms\  $n_0$.
Let $\ph \colon A \to C$ and $\rh \colon B \to C$
be unital \hm s, with $\rh$ surjective.
(As usual, we allow the case $C = 0$, $\ph = 0$, and $\rh = 0$.)
Then there exists a \rshd\  for $A \oplus_C B$ such that,
adopting the standard notation (from Definition~\ref{A2}),
the following hold:
\bit
\item[(1)]
The $l_0$-th stage algebra
(Definition~\ref{A2}) of this \rshd\  is $A$, and its obvious
\rshd\  of length $l_0$ is the given \rshd\  of $A$.
\item[(2)]
For $l_0 < k \leq l$, we have $\dim ( X_k) \leq \dim (X)$.
\item[(3)]
For $l_0 < k \leq l$, the matrix size $n (k)$ is equal to $n$.
\eit
In particular, $A \oplus_C B$ has a \rshd\  whose \tdim\  %
is $\max (d_0, \, \dim (X))$, and whose \mms\  is\linebreak
$\max (n_0, n)$.

Moreover, if $X^{(0)} \subset X$ is the closed subset of $X$ such that
\[
\Ker (\rh) =
 {\ts{ \left\{ f \in \Gm (E) \colon f |_{X^{(0)}} = 0 \right\} }},
\]
and if $X \setminus X^{(0)}$
is dense in $X$, then the \rshd\  may be chosen to have the additional
property:
\bit
\item[(4)]
For $l_0 < k \leq l$, we have $X_k \setminus X_k^{(0)}$ dense in $X_k$.
\eit
\end{prp}

\begin{pff}
We may clearly assume
$C = \Gm {\ts{ \left( \rsz{ E |_{X^{(0)}} } \right) }}$ with
$X^{(0)} \subset X$ closed, and that $\rh$ is the restriction map.
Using compactness of $X$ and local triviality, cover $X$ with
finitely many open sets $U_1, \dots, U_r$ such that
$E |_{\overline{U}_j}$ is trivial for $1 \leq j \leq r$.
We prove the proposition by induction on $r$.
For $r = 1$, we have $B \cong C (X, M_n)$, and
$A \oplus_C B$ has an obvious \rshd\  of length $l = l_0 + 1$.
Moreover, the properties~(1) through~(3) (and~(4), if applicable)
are immediate.

So suppose the proposition is known in all cases when $r$ sets
$U_1, \dots, U_r$ suffice, and suppose $E$ is given, with a
cover $\{ U_1, \dots, U_{r + 1} \}$ such that
$E |_{\overline{U}_j}$ is trivial for $1 \leq j \leq r + 1$.
Define
\[
\widetilde{X} = \bigcup_{j = 1}^r \overline{U}_j, \,\,\,\,\,\,
\widetilde{X}^{(0)} = \widetilde{X} \cap X^{(0)}, \,\,\,\,\,\,
\widetilde{B} = \Gm {\ts{ \left( \rsz{ E |_{\widetilde{X}} } \right) }},
\andeqn \widetilde{C} =
     \Gm {\ts{ \left( \rsz{ E |_{\widetilde{X}^{(0)}} } \right) }}.
\]
Let $\pi \colon C \to \widetilde{C}$ be the restriction map.
Let $\widetilde{\ph} = \pi \circ \ph \colon A \to \widetilde{C}$,
and let $\widetilde{\rh} \colon \widetilde{B} \to \widetilde{C}$
be given by restriction of sections from $\widetilde{X}$ to
$\widetilde{X}^{(0)}$.
Then $\widetilde{\ph}$ and $\widetilde{\rh}$ are unital, and
$\widetilde{\rh}$ is surjective.
By the induction hypothesis,
$A \oplus_{\widetilde{C}} \widetilde{B}$ has a \rsha\  such that
the properties~(1) through~(3) hold.

Define
\[
Y = \overline{X \setminus \widetilde{X}}
\andeqn
Y^{(0)} =
 Y \cap {\ts{ \left( \rsz{\widetilde{X} \cup X^{(0)} } \right) }}.
\]
Then $Y$ and $\widetilde{X}$ are closed subsets of $X$ which cover $X$.
\Wolog\  $Y \neq \varnothing$.
Now set
\[
B_0 = \Gm {\ts{ \left( \rsz{ E |_{Y} } \right) }}
  \andeqn
C_0 =
 \Gm {\ts{ \left( \rsz{ E |_{ Y^{(0)} } } \right) }}.
\]
Note that $B_0 \cong C (Y, M_n )$
because $Y \subset \overline{U}_{r + 1}$
and $E |_{\overline{U}_{r + 1}}$ is trivial.
Let $\rh_0 \colon B_0 \to C_0$ be the restriction map.
Define $\ph_0 \colon A \oplus_{\widetilde{C}} \widetilde{B} \to C_0$
as follows.
For $(a, f) \in A \oplus_{\widetilde{C}} \widetilde{B}$, we let
$\ph_0 (a, f)$ be the section in $C_0$ given by
\[
\ph_0 (a, f) |_{Y \cap \widetilde{X} } = f |_{Y \cap \widetilde{X} }
                      \andeqn
\ph_0 (a, f) |_{Y \cap  X^{(0)}} = \ph (a) |_{Y \cap  X^{(0)}}.
\]
The condition that $(a, f)$ be in
$A \oplus_{\widetilde{C}} \widetilde{B}$ ensures that we do in fact have
\[
f |_{Y \cap \widetilde{X} \cap X^{(0)} } =
 \ph (a)  |_{Y \cap \widetilde{X} \cap X^{(0)} }.
\]
We now claim that
\[
A \oplus_C B \cong
 {\ts{ \left( \rsz{ A \oplus_{\widetilde{C}} \widetilde{B} } \right) }}
       \oplus_{C_0} B_0.
\]
To see this, simply observe that the right hand side is the set of
all triples
\[
(a, f, g) \in A \oplus \widetilde{B}  \oplus B_0
\]
such that
\[
\ph (a)  |_{\widetilde{X}^{(0)}}  = f |_{\widetilde{X}^{(0)}},
           \,\,\,\,\,\,
\ph (a)  |_{Y \cap  X^{(0)}} = g |_{Y \cap  X^{(0)}},
                      \andeqn
f |_{Y \cap \widetilde{X} } = g |_{Y \cap \widetilde{X} }.
\]
The map that sends $(a, h) \in A \oplus_C B$ to
${\ts{ \left( \rsz{ a, h |_{\widetilde{X}}, \, h |_Y } \right) }}$
is now obviously an isomorphism.

The \ca\  %
${\ts{ \left( \rsz{ A \oplus_{\widetilde{C}} \widetilde{B} } \right) }}
       \oplus_{C_0} B_0$
has an obvious \rshd, with length $1$ larger than that for
$A \oplus_{\widetilde{C}} \widetilde{B}$.
Moreover, the conditions~(1) through~(3) clearly hold
(using Proposition~3.1.5 of \cite{Pr} for the dimension estimate
in~(2)).

It remains to show that if $X \setminus X^{(0)}$ is dense in $X$,
then also condition~(4) also holds.
By induction, we need only show that
$\widetilde{X} \setminus \widetilde{X}^{(0)}$ is dense in
$\widetilde{X}$ and $Y \setminus Y^{(0)}$ is dense in $Y$.
For the first of these, it suffices to show that each
$U_j \cap \left( \rsz{ X \setminus X^{(0)} } \right)$
is dense in $U_j$, since this set is then necessarily dense in
$\overline{U}_j$.
This is easy because $U_j$ is open and $X \setminus X^{(0)}$ is
dense in $X$.
For the second, it suffices to prove that
$\left( \rsz{ X \setminus \widetilde{X} } \right) \setminus Y^{(0)}$
is dense in $X \setminus \widetilde{X}$.
This is true because
\[
\ts{ \left( \rsz{ X \setminus \widetilde{X} } \right) }
    \setminus Y^{(0)}
  = \ts{ \left( \rsz{ X \setminus \widetilde{X} } \right) }
    \cap \ts{ \left( \rsz{ X \setminus X^{(0)} } \right) },
\]
and $X \setminus \widetilde{X}$ is open and $X \setminus X^{(0)}$ is
dense.
\end{pff}

\begin{cor}\label{A5.05}
Let $X$ be a \cpt\  space, let $E$ be a locally trivial continuous
field over $X$ with fiber $M_n$, and let $A = \Gm (E)$ be the
corresponding section algebra.
Then $A$ is a \rsha.
\end{cor}

In Section~3, we prove, in the separable case and
using more machinery, that $A \oplus_C B$ is a \rsha\  when
$B$ is merely required to be a \rsha.

To obtain the other
interesting corollaries, the following lemma is useful.

\begin{lem}\label{A5.1}
Let $A$, $B_1$, $B_2$, and $C$ be unital \ca s (allowing, however,
$C = 0$), and let $\ph \colon A \to C$ and $\rh \colon B_1 \oplus B_2 \to C$
be unital \hm s, with $\rh$ surjective.
Then there are unital \ca s $C_1$ and $C_2$ (possibly the zero \ca)
and unital \hm s
\[
\ph_1 \colon A \to C_1, \,\,\,\,\,\, \rh_1 \colon B_1 \to C_1, \,\,\,\,\,\,
\ph_2 \colon A \oplus_{C_1} B_1 \to C_2, \andeqn \rh_2 \colon B_2 \to C_2,
\]
with $\rh_1$ and $\rh_2$ surjective, such that
$A \oplus_C (B_1 \oplus B_2)
       \cong \left(A \oplus_{C_1} B_1 \right) \oplus_{C_2} B_2$.
\end{lem}

\begin{pff}
We have $C = C_1 \oplus C_2$, where $C_j = \rh (B_j)$.
Let $\rh_j \colon B_j \to C_j$ be the corresponding restriction of $\rh$.
Let $\pi_j \colon C \to C_j$ be the \pj\  map, and define
$\ph_1 = \pi_1 \circ \ph$ and $\ph_2 (a, b) = \pi_2 \circ \ph (a)$.
It is easy to check that the map $(a, b_1, b_2) \mapsto (a, b_1, b_2)$
is an isomorphism of the two different fibered products.
\end{pff}

\begin{cor}\label{A5.2}
In Definition~\ref{A1} we may replace $C (X, M_n)$ by an arbitrary
locally trivial homogeneous algebra without changing the class
of algebras defined.
Moreover, given a decomposition of this weaker
kind, there is a true decomposition with the same
\tdim\  and set of matrix sizes.
\end{cor}

\begin{pff}
We use Lemma~\ref{A5.1} to reduce to the case of constant rank,
which is then
covered by induction and Proposition~\ref{A5}.
\end{pff}

\begin{cor}\label{A5.3}
Let $A$ be a \rsha, and let $p \in A$ be a \pj.
Then $pAp$ is a \rsha,
with \tdim\  and \mms\  no larger than for $A$.
\end{cor}

\begin{pff}
Let $A$ have a decomposition in the standard notation
(as in Definition~\ref{A2}).
Let $p_k$ be the image of $p$ in $C_k$,
and let $p_k^{(0)}$ be the image of $p_k$ in $C_k^{(0)}$.
Then
\[
pAp = \left[ \cdots \rule{0em}{3ex} \left[ \left[
 p_0 C_0 p_0 \oplus_{p_1^{(0)} C_1^{(0)} p_1^{(0)} } p_1 C_1 p_1 \right]
 \oplus_{p_2^{(0)} C_2^{(0)} p_2^{(0)} } p_2 C_2 p_2 \right]
                      \cdots \right]
            \oplus_{p_l^{(0)} C_l^{(0)} p_l^{(0)} } p_l C_l p_l.
\]
Omitting terms for which $p_k = 0$,
the conclusion follows from Corollary~\ref{A5.2}.
\end{pff}

\begin{lem}\label{A7}
Let $A$ be a \rsha.
Then $M_n (A)$ is a \rsha.
It has a decomposition
with the same length and base spaces as for $A$, and
with \mms\  equal to $n$ times the \mms\  of $A$.
\end{lem}

\begin{pff}
This is immediate from the definition.
\end{pff}

\begin{cor}\label{A8}
Let $A$ be a \rsha.
Let $B$ be a unital \ca\  which is Morita equivalent to $A$.
Then $B$ is a \rsha.
\end{cor}

\begin{pff}
Proposition~2.1 of \cite{Rf0} implies that $B \cong p M_n (A) p$
for some $n$ and some \pj\  $p \in M_n (A)$.
So $B$ is a \rsha\  by Lemma~\ref{A7} and Corollary~\ref{A5.3}.
\end{pff}

\begin{lem}\label{A8.5}
Let $A$ be a \rsha\  with base spaces $X_k$, and let $Y$ be
a \cpt\  space.
Then $C (Y) \otimes A$ is a \rsha, and has a decomposition
with the same length as for $A$ and base spaces $Y \times X_k$.
In particular, the \tdim\  of this decomposition is at most
$\dim (Y)$ times the \tdim\  of $A$, and its \mms\  is the same
as for $A$.
\end{lem}

\begin{pff}
The first part is immediate from the definition.
(It is also a special case of Theorem~3.4 of \cite{Pd}.)
For the last sentence one needs
$\dim (Y \times X_k) \leq \dim (Y) + \dim (X_k)$, which is
Proposition~3.2.6 of \cite{Pr}.
\end{pff}

\section{Characterization of \rsha s}

In this section, we show that a separable unital \ca\  has a
\rshd\  with finite topological dimension \ifo\  there is a
finite upper bound on the dimensions of irreducible representations
and, for each $n$, the space of irreducible representations of
dimension $n$ is finite dimensional.
This result shows that there are many \rsha s.
We further prove a related (but not as nice)
characterization of \rsha s in the absence
of finite dimensionality.

In this section, we use the {\emph{modified}} covering dimension
for topological spaces.
This is the same as covering dimension, except that one only
considers covers consisting of complements of the zero sets of
continuous real valued functions on $X$.
See Definition~10.1.3 of \cite{Pr}, but note that this function
is denoted $\partial{\mathrm{im}} (X)$ there, and we will simply
write $\dim (X)$.
For compact Hausdorff spaces, and more generally for normal spaces,
this dimension agrees with the usual covering dimension
(Proposition~10.1.6 of \cite{Pr}), and for completely regular
spaces this definition gives $\dim (X) = \dim (\bt X)$
(Theorem~10.1.4 of \cite{Pr}).
The difference
will be relevant only in some of the results which are stated for
not necessarily $\sm$-compact spaces.

\begin{lem}\label{A9}
Let $A$ be a \rsha\  with total space $X$.
Then the map $x \mapsto \Ker (\ev_x)$ defines
a \ct\  bijection
\[
\coprod_{k = 0}^l
  \ts{ \left( \rsz{ X_k \setminus X_k^{(0)} } \right) } \to \Prim (A)
\]
whose restriction to each $X_k \setminus X_k^{(0)}$ is
a homeomorphism onto its image.
(We take $X_0^{(0)} = \varnothing$.)
In particular, every irreducible representation of $A$ is unitarily
equivalent
to a representation of the form $\ev_x$ for some $x \in X$.
\end{lem}

\begin{pff}
The proof is by induction on the length $l$.
If the length is zero, then $A = C(X, M_n)$ for some $X$ and $n$, and
the result is obvious.
So assume that the result is known for length $l$.
Let $A = B \oplus_{C ( X^{(0)}, \, M_n )} C (X, M_n)$,
where $B$ is a \rsha\  of length $l$.
Consider the exact sequence
\[
0 \longrightarrow
    C_0 {\ts{ \left( \rsz{ X \setminus X^{(0)}, \, M_n } \right) }}
   \longrightarrow A \longrightarrow B \longrightarrow 0,
\]
in which the first map is $f \mapsto (0, f)$ and the second
is $(a, f) \mapsto a$.
It shows that $\Prim (A)$ is the (nontopological)
disjoint union of the closed set $\Prim (B)$
and the open set
$\Prim {\ts{ \left( C_0 \left( \rsz{ X \setminus X^{(0)} },
           \, M_n \right)  \right) }} \cong X \setminus X^{(0)}$,
and these two identifications are homeomorphisms onto
their images.
So the result for $B$ implies the result for $A$.
\end{pff}

\medskip

We warn that, as one can see from the examples, not every
$\ev_x$ is irreducible.

The next definition is a generalization of Definition~3.5.7 of \cite{Hs}
(which is only stated for vector bundles), and is
presumably known in the generality in which we give it.
I am grateful to Larry Brown for the suggestion that the finite type
condition might be relevant here.

\begin{dfn}\label{N1}
Let $X$ be a topological space and let $E$ be a
locally trivial fiber bundle over $X$.
Then $E$ is said to have {\emph{finite type}} of there is
a finite open cover $U_1, \, U_2, \dots , U_n$ of $X$
such that $E |_{U_k}$ is trivial for $1 \leq k \leq n$.
\end{dfn}

We prove that vector bundles over sufficiently nice
finite dimensional spaces have finite type.
This is proved in Lemma~2.19 of \cite{Ml}, but that
reference is possibly not widely available and in any
case uses a slightly different definition of covering dimension.

\begin{lem}\label{N2}
Let $X$ be a compact Hausdorff space with $\dim (X) \leq d < \infty$.
Let $E$ be a vector bundle over $X$ with rank $r$.
Then there exist $d + 1$ open sets
$V_1, \, V_2, \dots , V_{d + 1} \subset X$
which cover $X$ and
such that $E |_{V_k}$ is trivial for $1 \leq k \leq d + 1$.
\end{lem}

\begin{pff}
The method of proof of Theorem~2.5 of \cite{Gd}
shows that there is a finite complex $Y$
with $\dim (Y) \leq d$, a \ct\  function $f \colon X \to Y$,
and a vector bundle $F$ over $Y$ with rank $r$, such
that $f^* (F) \cong E$.
(The argument reduces to the fact that $X$ is an inverse
limit of compact metric spaces of dimension at most $d$,
and such a space is in turn an inverse
limit of finite complexes of dimension at most $d$.)
Now Theorem~1.2.6 of \cite{Hs} provides $d + 1$ open sets
$W_1, \, W_2, \dots , W_{d + 1} \subset Y$
which cover $Y$ and
such that each connected component of each $W_k$ is
contractible.
Therefore $F |_{W_k}$ is trivial.
The proof is completed by taking $V_k = f^{-1} (W_k)$ and noting
that $E |_{V_k} \cong f^* \left( F |_{W_k} \right)$.
\end{pff}

\begin{lem}\label{N3}
Let $X$ be a locally compact $\sm$-compact Hausdorff space.
Then there exist open sets $U_1, \, U_2, \dots \subset X$
which cover $X$, whose closures are compact, and such that
$\overline{U}_j \cap \overline{U}_k = \varnothing$ for $| j - k | > 1$.
\end{lem}

\begin{pff}
Let $K_1, \, K_2, \dots \subset X$ be a sequence of
compact sets which covers $X$.
We construct $U_n$ by induction on $n$, satisfying:
\bit
\item[(1)]
$\bigcup_{k = 1}^n U_k \supset K_n$.
\item[(2)]
$\partial \left( \bigcup_{k = 1}^{n - 1} U_k \right) \subset U_n$.
\item[(3)]
$\overline{U}_j \cap \overline{U}_k = \varnothing$
for $| j - k | > 1$ and $1 \leq j, \, k \leq n$.
\item[(4)]
$\overline{U}_k$ is compact for $1 \leq k \leq n$.
\eit
We start the induction by choosing $U_1$ to be any open set which
contains $K_1$ and has compact closure.

Suppose given $U_1, \, U_2, \dots , U_n$ as above.
Choose an open set $V$ which
contains $K_{n + 1} \cup \bigcup_{k = 1}^{n} \overline{U}_k$
and has compact closure.
Define
\[
W = V \setminus \bigcup_{k = 1}^{n - 1} \overline{U}_k
  \andeqn
L = \left( K_{n + 1} \setminus \bigcup_{k = 1}^{n} U_k  \right)
  \cup \partial \left( \bigcup_{k = 1}^{n} U_k \right).
\]
Thus $W$ is open and $L$ is compact.
Also $L \subset W$, because
\[
V \setminus \bigcup_{k = 1}^{n} U_k
 \subset
V \setminus \bigcup_{k = 1}^{n - 1} \overline{U}_k
                 \andeqn
\partial \left( \bigcup_{k = 1}^{n} U_k \right)
 \subset
\bigcup_{k = 1}^{n} \overline{U}_k \setminus \bigcup_{k = 1}^{n} U_k
 \subset W.
\]

Now choose an open set $U_{n + 1}$ with compact closure such that
$L \subset U_{n + 1} \subset \overline{U}_{n + 1} \subset W$.
Conditions~(1), (2), and~(4) are satisfied for $n + 1$ by
construction.
For~(3), we need only show that
$\overline{U}_{n + 1} \cap \overline{U}_k = \varnothing$
for $1 \leq k \leq n - 1$.
But this follows from the relation
$W \cap \bigcup_{k = 1}^{n - 1} \overline{U}_k = \varnothing$.
This completes the induction.
The lemma is now proved because the sets $K_n$ cover $X$.
\end{pff}

\begin{lem}\label{N4}
Let $X$ be a locally compact $\sm$-compact Hausdorff space
with\linebreak
$\dim (X) \leq d < \infty$.
Let $E$ be a vector bundle over $X$ with (constant) rank $r$.
Then $E$ has finite type.
\end{lem}

\begin{pff}
Choose open sets $U_1, \, U_2, \dots \subset X$ satisfying the
conclusion of Lemma~\ref{N3}.
We have $\dim \left( \overline{U}_k \right) \leq d$
by Proposition~3.1.5 of \cite{Pr}.
(Since $X$ and $\overline{U}_k$ are normal, the modified
covering dimension is the same as the usual one.)
By Lemma~\ref{N2} there exist open subsets
$V_{k, 1}, \, V_{k, 2}, \dots , V_{k, d + 1} \subset \overline{U}_k$
which cover $\overline{U}_k$ and
such that $E |_{V_{k, l}}$ is trivial for $1 \leq l \leq d + 1$.
Now define
\[
W_l^{(0)} = \bigcup_{k \,\, {\text{even}}} (U_k \cap V_{k, l})
  \andeqn
W_l^{(1)} = \bigcup_{k \,\, {\text{odd}}} (U_k \cap V_{k, l}).
\]
Each of these $2 (d + 1)$ sets is the union of disjoint
open subsets of $X$ over which $E$ is trivial,
since $U_j \cap U_k = \varnothing$ for $|j - k| > 1$.
Therefore $E$ is trivial over each $W_l^{(i)}$.
Since these sets cover $X$, the proof is complete.
\end{pff}

\medskip

The following equivalent conditions for finite type will be very useful.

\begin{lem}\label{N5}
Let $X$ be a completely regular space, and
let $E$ be a vector bundle over $X$ with (constant) rank $r$.
Then \tfae:
\bit
\item[(1)]
$E$ has finite type.
\item[(2)]
There exists a vector bundle $H$ over $X$ such that
$E \oplus H$ is trivial.
\item[(3)]
There exists a vector bundle $F$ over $\bt X$,
the Stone-\v{C}ech compactification of $X$,
such that $F |_X \cong E$.
\item[(4)]
There exists a compactification $Y$ of $X$,
and a vector bundle $F$ over $Y$,
such that $F |_X \cong E$.
\eit
\end{lem}

\begin{pff}
The equivalence of~(1) and~(2) is the equivalence of
conditions~(1) and~(3) of Proposition~3.5.8 of \cite{Hs}.
That~(3) implies~(4) is trivial.
For~(4) implies~(2), observe that $F$ is a direct summand
of a trivial bundle
because $Y$ is compact.
Restricting to $X$, we see that $E$
is a direct summand of a trivial bundle.

It remains to show that~(2) implies~(3).
We may assume $E \oplus H = X \times \C^n$ for some $n$.
Equip $X \times \C^n$ with the usual scalar product,
and let $p$ be the orthogonal \pj\  from $X \times \C^n$ to $E$.
We may think of $p$ as a bounded \ct\  function
from $X$ to $M_n$.
Then there is a \ct\  function
$q$ from $\bt X$ to $M_n$ which extends $p$.
The values of $q$ must be \pj s, so $F = q (X \times \C^n)$
is a vector bundle over $\bt X$ such that
$F |_X \cong p (X \times \C^n) \cong E$.
\end{pff}

\medskip

Our next goal is a similar characterization of finite type for
locally trivial $M_n$-bundles.
We will in fact see that an $M_n$-bundle has finite type \ifo\  the
underlying $n^2$-dimensional vector bundle does.
We need several lemmas.

\begin{lem}\label{N6}
Let $X$ be a compact Hausdorff space.
Let $E$ be a vector bundle over $X$ with rank $r > 0$.
Then there exists a vector bundle $F$ over $X$ with rank $s > 0$
such that $E \otimes F$ is trivial.
\end{lem}

\begin{pff}
Use Theorem~10.1 in Chapter~10 of \cite{ES}
to write $X = \invlim X_{\af}$ for finite complexes $X_{\af}$.
Then there is $\af$ and a vector bundle $E_0$ over $X_{\af}$
such that, with $f \colon X \to X_{\af}$ being the canonical map, we have
$E \cong f^* (E_0)$.
Write $X_{\af}$ as a disjoint union $X_{\af} = \coprod_{k = 1}^n Y_k$
of connected finite complexes.
Lemma~12 of \cite{DK} provides vector bundles $G_k$ over $Y_k$
with positive rank $s (k)$ such that $G_k \otimes E |_{Y_k}$ is trivial.
Set $t (k) = \prod_{l \neq k} s (l)$.
Then
$F_0 = \coprod_{k = 1}^n G_k \otimes
        {\ts{ \left( \rsz{ Y_k \times \C^{t (k)} } \right) }}$
is a vector bundle over $X_{\af}$
with constant rank $s = \prod_{k = 1}^n s (k) > 0$ and
such that $E_0 \otimes F_0$ is trivial.
Take $F = f^* (F_0)$.
\end{pff}

\begin{lem}\label{N7}
Let $X$ be a compact Hausdorff space, let $S$ be a dense
subset of $X$, and let $n, \, N \in \N$.
Let $x \mapsto P_x$ be a \ct\  function from $X$ to the
space $L (M_N)$ of bounded operators on $M_N$.
Suppose that, for all $x \in S$, the map $P_x$ is a
conditional expectation from $M_N$ onto
a subalgebra $A_x \cong M_n$ of $M_N$.
Then $P_x$ is a
conditional expectation from $M_N$ onto
a subalgebra $A_x \cong M_n$ for every $x \in X$.
Moreover
the assignment $x \mapsto A_x$ defines a locally
trivial \ct\  field of \ca s over $X$ with fiber $M_n$.
\end{lem}

\begin{pff}
It is easy to check that the set of $x \in X$ for which
the range of $P_x$ is not a subalgebra is open in $X$,
and hence empty.
The other algebraic and norm properties of
conditional expectations extend over $X$ in the same
way.
Similarly, the trace of $P_x$ is $n^2$ for all $x \in S$,
hence for all $x \in X$.
It follows that $P_x$ has rank exactly equal to $n^2$ for
all $x \in X$.

Let
\[
V = \{ x \in X \colon A_x \not\cong M_n \} \subset X,
\]
and let $x \in V$.
Since $\dim (A_x) = n^2$, it follows that there are
$n + 1$ nonzero \mops\  $p_1, \dots, p_{n + 1} \in A_x$.
Define $a_k (y) = P_y (p_k)$ for $y \in X$.
Then the $a_k$ are continuous functions from $X$ to $M_N$,
such that $a_k (y) \in A_y$ for all $y$.
Standard functional calculus techniques (or Chapter~14 of \cite{Lr})
show that, on some \nbhd\  of $x$, there are
\ct\   \mops\  %
$y \mapsto q_1 (y), \dots, y \mapsto q_{n + 1} (y)$,
whose values at $y$ are in the C*-subalgebra of $M_N$
generated by $a_1 (y), \dots, a_{n + 1} (y)$.
(This is semiprojectivity of $\C^{n + 1}$.)
Therefore, for $y$ close enough to $x$,
there are $n + 1$
nonzero \mops\  $q_1 (y), \dots, q_{n + 1} (y) \in A_y$.
It follows that $A_y \not\cong M_n$, so that $V$ is open.
Since $V \cap S = \varnothing$, we have
$A_x \cong M_n$ for all $x \in X$.

It remains to prove local triviality.
Let $x \in X$.
Choose a system $(e_{j, k})_{1 \leq j, \, k \leq n}$ of
matrix units for $A_x$.
Using the conditional expectation $P_x$ and
standard functional calculus techniques (or \cite{Lr}) as in
the previous paragraph, we may find a \nbhd\  $U$ of $x$
and \ct\  functions $f_{j, k} \colon V \to M_N$ for
$1 \leq j, \, k \leq n$ such that
$f_{j, k} (y) \in A_y$, $f_{j, k} (x) = e_{j, k}$, and
$(f_{j, k} (y))_{1 \leq j, \, k \leq n}$ is a system of
matrix units.
(This is semiprojectivity of $M_n$.)
Reducing the size of $V$ if necessary, we may find a
\ct\  function $s \colon V \to M_N$ such that
\[
s (y)^* s (y) = e_{1, 1}, \,\,\,\,\,\,
s (y) s (y)^* = f_{1, 1} (y),  \andeqn s (x) = e_{1, 1}.
\]
Then $y \mapsto w (y) = \sum_{j = 1}^n f_{j, 1} (y) s (y) e_{1, j}$
is a continuous function from $V$ to the unitaries in $M_N$
such that $w e_{j, k} w^* = f_{j, k} (y)$.
Thus $w A_x w^* = A_y$, and the bundle $y \mapsto A_y$
is trivial over $V$.
\end{pff}

\begin{prp}\label{N8}
Let $X$ be a completely regular space, and
let $A$ be a locally
trivial \ct\  field of \ca s over $X$ with fiber $M_n$.
Then \tfae:
\bit
\item[(1)]
$A$ has finite type.
\item[(2)]
$A$ has finite type when regarded as a complex vector
bundle over $X$ by forgetting structure.
\item[(3)]
There exists a locally
trivial \ct\  field $B$ of \ca s over $\bt X$
such that $B |_X \cong A$.
\item[(4)]
There exists a compactification $Y$ of $X$,
and a locally
trivial \ct\  field $B$ of \ca s $F$ over $Y$,
such that $B |_X \cong A$.
\item[(5)]
There exists $k$ and a locally
trivial \ct\  field $C$ of \ca s over $X$ with fiber
$M_k$ such that $A \otimes C$ is trivial.
\eit
\end{prp}

\begin{pff}
We prove
$(1) \Rightarrow (2) \Rightarrow (5) \Rightarrow (3) \Rightarrow (4)
       \Rightarrow (1)$.
The implications $(1) \Rightarrow (2)$ and $(3) \Rightarrow (4)$
are immediate, and $(4) \Rightarrow (1)$ follows directly
from compactness of $Y$ and local triviality of $B$.

We prove $(2) \Rightarrow (5)$.
Define a scalar product in the fibers of $A$ by taking
the Hilbert-Schmidt scalar product.
Since $A$ has finite type as a vector bundle,
$(1) \Rightarrow (4)$ of Lemma~\ref{N5} provides a compactification
$Y$ of $X$ and a vector bundle $F$ over $Y$ such that $F |_X \cong A$.
Lemma~\ref{N6} then provides a vector bundle $H$ over $Y$
and an isomorphism $a \colon F \otimes H \to Y \times \C^N$ for some $N$.
We may assume $H$ has a continuously varying scalar
product and equip $Y \times \C^N$ with the usual scalar product.
Then, replacing $a$ by $a (a^* a)^{-1/2}$, we may
assume that $a$ is unitary.

Let $A^{\op}$ denote the locally
trivial \ct\  field over $X$ with fiber
$M_n$ obtained by reversing the multiplication in every fiber.
Then, as is well known,
there is a *-isomorphism $A \otimes A^{\op} \cong L (A)$,
the locally trivial \ct\  field over $X$ with fiber
$M_{n^2} \cong L (A_x)$ using the Hilbert space structure
on $A_x$ from the previous paragraph.
(Lacking a reference, we give a formula.
Define $\pi_x \colon A_x \otimes A_x^{\op} \to L (A_x)$ by
$\pi_x (a_1 \otimes a_2) (c) = a_1 c a_2$.
One easily checks that this defines a unital *-\hm.
Simplicity of the domain $A_x \otimes A_x^{\op}$ and a dimension
count show that it is bijective.)
With similar notation, we have
\[
L (A) \otimes L (H |_X) \cong
  L {\ts{ \left( \rsz{ A \otimes H |_X } \right) }} \cong
  L {\ts{ \left( \rsz{ X \times \C^N } \right) }} \cong X \times M_N.
\]
This is~(5) with $C = A^{\op} \otimes L (H |_X)$.

It remains to prove $(5) \Rightarrow (3)$.
By hypothesis, there exist $k$ and a
locally trivial \ct\  field $C$ over $X$ with fiber
$M_k$ such that $A \otimes C \cong M_{nk}$.
For each $x \in X$, the tensor product decomposition
$M_{nk} \cong A_x \otimes C_x$ gives a
conditional expectation $P_x \colon M_{nk} \to A_x$.
Then $x \mapsto P_x$ is a \ct\  function from
$X$ to the operators of norm $1$ on $M_{nk}$.
Therefore it extends to a \ct\  function,
still called $x \mapsto P_x$, defined on $\bt X$.
By Lemma~\ref{N7}, the assignment $x \mapsto P_x (M_{nk})$
defines a locally trivial \ct\  field over $\bt X$ with fiber
$M_n$, whose restriction to $X$ is obviously $A$.
\end{pff}

\medskip

Several more lemmas are needed before we can prove the
characterization theorems.

\begin{lem}\label{N9}
Let $X$ be a second countable locally compact Hausdorff space,
and let $Y$ be a compactification of $X$.
Let $A$ be a locally
trivial \ct\  field of \ca s over $Y$ with fiber $M_n$,
and let $D$ be a separable C*-subalgebra of the section
algebra $\Gm \left( A |_{Y \setminus X} \right)$.
Then there are a compactification $Z$ of $X$,
a surjective map $h \colon Y \to Z$ which is the identity on $X$,
and a locally trivial \ct\  field $B$ over $Z$ with fiber $M_n$,
such that:
\bit
\item[(1)]
$Z$ is second countable.
\item[(2)]
$\dim (Z) = \dim (X)$.
\item[(3)]
There is an isomorphism
$\ph \colon h^* (B) |_{Y \setminus X} \cong A |_{Y \setminus X}$
such that the  range of the induced map
$\sm \colon \Gm \left( B |_{Z \setminus X} \right)
              \to \Gm \left( A |_{Y \setminus X} \right)$
contains $D$.
(Here $\sm (b) (x) = \ph_x (b ( h (x)))$.)
\eit
\end{lem}

\begin{pff}
We first find a compactification $Z_0$ with $h_0 \colon Y \to Z_0$
such that all of the conclusion except part~(2) holds.
By~(4) implies~(5) of Lemma~\ref{N8}, applied to the compact space $Y$,
there exists a locally trivial \ct\  field $C$ over $Y$ with fiber $M_k$
and an isomorphism $\mu_0 \colon A \otimes C \to Y \times M_{nk}$.
This map induces a \hm\  %
$\mu \colon \Gm \left( A |_{Y \setminus X} \right)
                         \to C (Y \setminus X, \, M_{nk})$.
It also gives a \ct\  function $x \mapsto P_x$ from $Y$ to the set of
conditional expectations from $M_{nk}$ onto subalgebras, such that
$P_x ( M_{nk}) = (\mu_0)_x (A_x \otimes 1)$ for all $x \in Y$.
For $x, \, y \in Y \setminus X$ define $x \sim y$ to mean
$P_x = P_y$ and $\mu (a) (x) = \mu (a) (y)$ for all $a \in D$.
This defines an equivalence relation on $Y \setminus X$, and we extend
it to all of $Y$ by taking $x \sim y$ exactly when $x = y$ for
$x, \, y \in X$.
Let $Z_0$ be the maximal ideal space of the commutative unital \ca\  %
\[
R = \{ f \in C (Y) \colon f (x) = f (y) \,\, {\text{whenever}} \,\, x \sim y \}.
\]
Let $h_0 \colon Y \to Z_0$ be the \ct\  surjective map induced by the
inclusion of $R$ in $C (Y)$.
Note that $h_0$ sends $X$ homeomorphically onto $h_0 (X) \subset Z_0$,
and we can thus identify $X$ with $h_0 (X)$.
Moreover, $h_0 (x) = h_0 (y)$ \ifo\   $x \sim y$.

The map $x \mapsto P_x$ defines a \ct\  function
$x \mapsto \overline{P}_x$ from $Z_0$ to the set of
conditional expectations from $M_{nk}$ onto subalgebras.
It follows from Lemma~\ref{N7} that $x \mapsto \overline{P}_x ( M_{nk})$
defines a locally trivial \ct\  field $B_0$ over $Z_0$ with fiber $M_n$.
Since $x \mapsto P_x ( M_{nk})$ is just $A$
(the isomorphism being obtained from $\mu_0$),
there is an obvious isomorphism $h_0^* (B_0) \cong A$.
The range of the induced map from
$\Gm \left( B_0 |_{Z_0 \setminus X}\right)$
to $\Gm \left( A |_{Y \setminus X}\right)$
consists exactly of all sections $a$
of $A |_{Y \setminus X}$ such that
$(\mu_0)_x (a (x)) = (\mu_0)_y (a (y))$ whenever $x \sim y$.
In particular, it contains $D$, as required.

We show $Z_0$ is second countable.
Extend $x \mapsto P_x$ to a continuous function on $Y$
(whose values need no longer be conditional expectations).
Let
\[
\widetilde{D}
   = \{ a \in C (Y, \, M_{nk}) \colon a |_{Y \setminus X} \in \mu (D) \}.
\]
Then $\widetilde{D}$ is separable because $X$ is second countable and
there is an exact sequence
\[
0 \longrightarrow C_0 (X, \, M_{nk}) \longrightarrow  \widetilde{D}
          \longrightarrow D \longrightarrow 0.
\]
The extended function $x \mapsto P_x$ and the elements of
$\widetilde{D}$ are all images of functions on $Z_0$,
and together they separate the points of $Z_0$.
This implies that $Z_0$ is second countable.

We now modify $Z_0$ so as to obtain condition~(2) as well.
By Proposition~10.3.11 of \cite{Pr}, there exist a
compact Hausdorff space $Z$ and continuous surjections
$h \colon Y \to Z$ and $g \colon Z \to Z_0$ such that $h_0 = g \circ h$,
$\dim (Z) = \dim (Y)$, and $Z$ is also second countable.
The proof is completed by taking $B = g^* (B_0)$.
\end{pff}

\begin{lem}\label{N10}
Let $X$ be a locally compact Hausdorff space, and
let $A$ be a locally trivial \ct\  field over $X$ with fiber $M_n$
which has finite type.
Let $\Gm_0 (A)$ be the \ca\  of continuous sections of $A$ which
vanish at infinity on $X$.
Let
\[
0 \longrightarrow \Gm_0 (A) \longrightarrow C \longrightarrow D
                \longrightarrow 0
\]
be an exact sequence with $C$ and $D$ unital.
Then there is a compactification $Y$ of $X$, a
locally trivial \ct\  field $B$ over $Y$ with fiber $M_n$,
and a unital \hm\  %
$\ph \colon D \to \Gm \left( B |_{Y \setminus X}\right)$ such that:
\bit
\item[(1)]
$\dim (Y) = \dim (X)$.
\item[(2)]
$B |_X \cong A$.
\item[(3)]
There is an isomorphism
$\ps \colon D \oplus_{\Gm ( B |_{Y \setminus X} )} \Gm (B) \to C$
which fits into a commutative diagram with exact rows
\[
\begin{array}{ccccccccc}
0 & \longrightarrow & \Gm_0 (B |_X) & \longrightarrow
      & D \oplus_{\Gm ( B |_{Y \setminus X} )} \Gm (B)
      & \longrightarrow & D & \longrightarrow & 0
                                 \\
& &  \,\,\,\,\,\, \downarrow \cong & & \,\,\,\,\,\, \downarrow  \ps & &
             \,\,\,\,\,\, \downarrow  = & &
                                 \\
0 & \longrightarrow & \Gm_0 (A) & \longrightarrow & C &
        \longrightarrow & D & \longrightarrow & 0
\end{array}
\]
in which the first vertical map is induced by the isomorphism
of~(2) and the pullback is via $\ph$ and the restriction of
sections to $Y \setminus X \subset X$.
\eit

Moreover, if $X$ is second countable and $D$ is separable, then
$Y$ may be chosen to be second countable.
\end{lem}

\begin{pff}
We begin by identifying the multiplier algebra of $\Gm_0 (A)$.
By parts~(3) and~(5) of Proposition~\ref{N8},
there is a locally trivial \ct\  field $B_0$
over $\bt X$ with fiber $M_n$ such that $B_0 |_X \cong A$,
and, moreover, $B_0$ is a subbundle of a trivial bundle
$\bt X \times M_N$ for some $N$.
All bounded \ct\  sections of $A$ extend uniquely to sections
of $\bt X \times M_N$, and clearly the sections one gets this way
are exactly the \ct\  sections of $B_0$.
It now follows from Theorem~3.3 of \cite{APT} that
$M (\Gm_0 (A) )$ can be identified with $\Gm (B_0)$.

The exact sequence of the hypothesis thus yields \hm s
\[
\ta \colon C \to M (\Gm_0 (A) ) = \Gm (B_0)  \andeqn
\overline{\ta} \colon
    D \to \Gm (B_0) / \Gm_0 (A)
       = \Gm {\ts{ \left( B_0 |_{ \bt X \setminus X} \right) }}.
\]
In the general case, we take $Y = \bt X$, $B = B_0$, and
$\ph = \overline{\ta}$.
We have $\dim (\bt X) = \dim (X)$ by Theorem~10.1.4 of \cite{Pr}
(recall our definition of $\dim$).
We now construct  $\ps$.
Let $\pi \colon C \to D$ be the quotient map, and let
$\rh \colon
 \Gm (B_0) \to \Gm {\ts{ \left( B_0 |_{ \bt X \setminus X} \right) }}$
be restriction.
By definition, if $c \in C$ then
$\rh \circ \ta (c) = \ph \circ \pi (c)$, so
$\ps = (\ph, \rh)$ is a well defined \hm\  making the
diagram in~(3) commute.
Moreover, $\ps$ is an isomorphism by the Five Lemma.

In the second countable case, we apply Lemma~\ref{N9} to $B_0$, $\bt X$,
and
$\overline{\ta} (D) \subset
     \Gm {\ts{ \left( B_0 |_{ \bt X \setminus X} \right) }}$,
obtaining a second countable
compactification $Y$ of $X$ with $\dim (Y) = \dim (\bt X) = \dim (X)$,
and a locally trivial \ct\  field $B$ over $Y$, such that
$\Gm {\ts{ \left( B |_{Y \setminus X} \right) }}$
can be canonically identified with a subalgebra of
$\Gm {\ts{ \left( B_0 |_{ \bt X \setminus X} \right) }}$
which contains $\overline{\ta} (D)$.
These changes from the case above do not change the pullback
$D \oplus_{\Gm ( B |_{Y \setminus X} )} \Gm (B)$, so the lemma is
proved in this case also.
\end{pff}

\medskip

We let $\Prim (A)$ denote the primitive ideal space of a \ca\  $A$,
and we let $\Prim_n (A)$ denote the subspace of $\Prim (A)$ consisting
of the kernels of $n$-dimensional representations of $A$.
We summarize some of the standard facts.
We refer to $E$ as in~(4) as the induced \ct\  field over
$\Prim_n (A)$.

\begin{thm}\label{N10.5}
For each finite $n$:
\bit
\item[(1)]
$\bigcup_{k \leq n} \Prim_k (A)$ is closed in $\Prim (A)$.
\item[(2)]
$\Prim_n (A)$ is open in $\bigcup_{k \leq n} \Prim_k (A)$.
\item[(3)]
$\Prim_n (A)$ is locally compact Hausdorff.
\item[(4)]
There is
a locally trivial \ct\  field $E$ over $\Prim_n (A)$ with fiber $M_n$
such that the subquotient of $A$ corresponding to $\Prim_n (A)$ is
isomorphic to $\Gm_0 (E)$.
\eit
\end{thm}

\begin{pff}
Parts~(1) and~(2) are Proposition 3.6.3~(i) of \cite{Dx},
and part~(3) is Proposition 3.6.4~(i) of \cite{Dx}.
Part~(4) is Theorem~3.2 on page~249 of \cite{Fl}.
(Also see Theorems~3 and~5 of \cite{TT}.)
\end{pff}

\begin{prp}\label{N11}
Let $A$ be a unital \ca.
Suppose that there is $N \in \N$ such that all
irreducible representations of $A$ have dimension at most $N$,
and suppose that the induced \ct\  fields over the subspaces
$\Prim_n (A)$ (Theorem \ref{N10.5}~(4)) all have finite type.
Then $A$ has a \rshd\  which, when given in the standard notation
(see Definition~\ref{A2}) has the following properties:
\bit
\item[(1)]
$n (0) \leq n (1) \leq \cdots \leq n (l)$.
\item[(2)]
$\dim (X_k) \leq \dim \left( \Prim_{n (k)} (A) \right)$
for $0 \leq k \leq l$.
\item[(3)]
$X_k \setminus X_k^{(0)}$ is dense in $X_k$ for $1 \leq k \leq l$.
\eit
In particular, it has \tdim\  at most
$\max_{1 \leq n \leq N} \dim ( \Prim_n (A) )$ and \mms\  at most $N$.

If $A$ is separable, then the \rshd\  can be chosen so that,
in addition:
\bit
\item[(4)]
Every \ca\  in the \rshd\  is separable.
\eit
\end{prp}

\begin{pff}
Using Proposition~\ref{A5} and Theorem~\ref{N10.5},
this follows from Lemma~\ref{N10} by induction.
(The application of part~(4) of the conclusion of Proposition~\ref{A5}
is justified as follows:
In the application of Lemma~\ref{N10}, we use a compactification
of $\Prim_n (A)$, and $\Prim_n (A)$ is dense in a compactification
by definition.)
\end{pff}

\begin{lem}\label{N12}
Let $X$ be a metric space,
and let $E$ be a vector bundle over $X$ with rank $r$.
Suppose that there exist open sets $U_1, \dots, U_n \subset X$
and closed sets $L_1, \dots, L_n \subset X$, such that
the sets $U_k \cap L_k$ cover $X$ and
$E |_{U_k \cap L_k}$ has finite type.
Then $E$ has finite type.
\end{lem}

\begin{pff}
Applying the definition of finite type, we may immediately reduce
to the case in which
$E |_{U_k \cap L_k} \cong (U_k \cap L_k) \times \C^r$
(but possibly with $n$ larger).
Since $X$ is metric, $U_k$ is paracompact.
Moreover, $U_k \cap L_k$ is closed in $U_k$.
By paracompactness, $E |_{U_k}$ has a scalar product,
and by the usual methods we may assume that there is a
unitary isomorphism
$u \colon E |_{U_k \cap L_k} \to (U_k \cap L_k) \times \C^r$.
For each $x \in U_k$ choose an open set $V_x \subset U_k$ with
$x \in V_x$ and a unitary trivialization
$v_x \colon E |_{V_x} \to V_x \times \C^r$.
Then $v_x u^*$ defines a \ct\  function from $V_x \cap L_k$ to
the unitaries on $\C^r$, and so the Tietze Extension Theorem
gives a \ct\  function $a_x \colon V_x \to M_r$ which extends this
function.
Choose a locally finite partition of unity $( f_{\af} )_{\af \in I}$
subordinate to the open cover $( V_x )_{x \in U_k}$ of $U_k$,
with ${\mathrm{supp}} (f_{\af} ) \subset V_{x (\af)}$.
Define a morphism of vector bundles $a \colon E |_{U_k} \to U_k \times \C^r$
by $a = \sum_{\af \in I} f_{\af} a_{x (\af)}^* v_{x (\af)}$.
Then $a |_{U_k \cap L_k} = u$.
Let $W_k = \{ x \in U_k \colon a (x) \,\, {\text{is invertible}} \}$,
which is an open subset of $X$ containing $U_k \cap L_k$ and over
which $E$ is trivial.
Since the sets $W_k$ cover $X$, this shows that $E$ has finite type.
\end{pff}

\begin{thm}\label{N13}
Let $A$ be a separable unital \ca, and let $N \in \N$.
\Tfae:
\bit
\item[(1)]
$A$ has a \rshd\  with maximum matrix size at most $N$.
\item[(2)]
$A$ has a \rshd\  with maximum matrix size at most $N$
and whose total space is second countable.
\item[(3)]
All irreducible representations of $A$ have dimension at most
$N$, and for $1 \leq n \leq N$ the induced \ct\  field on
$\Prim_n (A)$ has finite type.
\eit
\end{thm}

\begin{thm}\label{N14}
Let $A$ be a separable unital \ca, and let $N, \, d \in \N$.
\Tfae:
\bit
\item[(1)]
$A$ has a \rshd\  with maximum matrix size at most $N$
and topological dimension at most $d$.
\item[(2)]
$A$ has a \rshd\  with maximum matrix size at most $N$ and
topological dimension at most $d$,
whose total space is second countable.
\item[(3)]
All irreducible representations of $A$ have dimension at most
$N$, and for $1 \leq n \leq N$ we have
$\dim (\Prim_n (A)) \leq d$.
\eit
\end{thm}

It is not entirely clear (beyond Proposition~\ref{N11}) what happens
in the absence of separability.

\medskip

{\emph{Proof of Theorems~\ref{N13} and~\ref{N14}:}}
In both cases, $(2) \Rightarrow (1)$ is trivial.
In Theorem~\ref{N13}, $(3) \Rightarrow (2)$ follows from
Proposition~\ref{N11}.
In Theorem~\ref{N14}, one uses in addition Lemma~\ref{N4}
and Proposition~\ref{N8}.

We prove $(1) \Rightarrow (3)$.
The condition on dimensions of representations is immediate
from Lemma~\ref{A9}.
For the other parts,
let a \rshd\  for $A$ be given as in the notation of Definition~\ref{A2}.
Apply Lemma~\ref{A9} and intersect everything with
$\Prim_n (A)$.
The result is a \ct\  bijection
\[
\coprod_{n (k) = n}
        {\ts{ \left( \rsz{ X_k \setminus X_k^{(0)} } \right) }}
   \to \Prim_n (A)
\]
whose restriction to each $X_k \setminus X_k^{(0)}$ is
a homeomorphism onto its image $Y_k$.
Moreover, $Y_k \subset \Prim_n (A)$ is locally compact, and hence
is the intersection of a closed set, say $L_k$, and an open
set, say $U_k$, in $\Prim_n (A)$.
The subquotient of $A$ corresponding to $Y_k$ is just
$C_0 {\ts{ \left( \rsz{ X_k \setminus X_k^{(0)}, \, M_n } \right) }}$.
With $E$ being the induced \ct\  field on $\Prim_n (A)$,
this shows that $E |_{Y_k}$ is trivial.
Now $\Prim_n (A)$ is second countable and locally compact (in
particular, regular), and so the version of the Urysohn Metrization
Theorem given in Theorem~4.4.1 of \cite{Mn} implies it is metrizable.
Therefore Lemma~\ref{N12} implies that $E$ has finite type.
This proves $(1) \Rightarrow (3)$ in Theorem~\ref{N13}.

For the corresponding part of Theorem~\ref{N14},
note that $X_k \setminus X_k^{(0)} \cong Y_k$ is an
$F_{\sm}$-set in $X_k$.
Therefore $\dim (Y_k) \leq \dim (X_k) \leq d$ by
the remark after Proposition~3.5.4 of \cite{Pr}.
Since $\Prim_n (A)$ is the union of finitely many such sets,
and each is also an $F_{\sm}$-set in $\Prim_n (A)$,
Proposition~3.5.3 of \cite{Pr} implies $\dim (\Prim_n (A)) \leq d$.
\QED

\section{Closure properties}

In the first section, we saw that the class of \rsha s is closed
under finite direct sums, tensor products with $M_n$ and commutative
\ca s, and passage to corners.
In this section we extend that list in the separable case, using
the characterization theorems of the previous section.
We show that quotients, pullbacks, and tensor products of
separable \rsha s are again separable \rsha.
Moreover, a quotient of a separable \rsha\  with
finite topological dimension is again an algebra of the same kind.
We also give an example of a separable unital \ca\  which has a
finite upper bound on the dimensions of its irreducible representations
but is not a \rsha.
This example is even a unital subalgebra of a separable homogeneous \ca.

The first several results are corollaries of
Theorems~\ref{N13} and~\ref{N14}.

\begin{prp}\label{N15}
Let $A$ be a separable \rsha, and let $I$ be an ideal in $A$.
Then $A/I$ is a \rsha.
Moreover, $A/I$ has a decomposition in which
the \tdim\  and \mms\  are no larger than for $A$.
\end{prp}

\begin{pff}
Condition~(3) in Theorems~\ref{N13} and~\ref{N14} passes to quotients.
(We have $\dim (\Prim_n (A/I)) \leq \dim (\Prim_n (A))$ by
Proposition~3.1.5 of \cite{Pr}, because $\Prim_n (A/I)$ is a closed
subset of $\Prim_n (A)$.)
\end{pff}

\begin{prp}\label{N16}
(Compare with Theorem~11.4 of \cite{Pd}.)
Let $A$ and $B$ be separable \rsha s, and let $\ph \colon A \to C$
and $\rh \colon B \to C$ be \hm s  with $\ph$ unital and $\rh$ surjective.
Then $A \oplus_C B$ is a \rsha.
Moreover, $A \oplus_C B$ has a decomposition in which
the \tdim\  and \mms\  are each no larger than the maximum of the
corresponding quantities for $A$ and $B$.
\end{prp}

\begin{pff}
We have an exact sequence
\[
0 \longrightarrow \Ker (\rh)
   \longrightarrow A \oplus_C B \longrightarrow B \longrightarrow 0
\]
(as in the proof of Lemma~\ref{A9}).
So $\Prim (A \oplus_C B)$ is the (nontopological)
disjoint union of the open set
$\Prim ( \Ker (\rh) )$ (also an open subset of $\Prim (A)$)
and the closed set $\Prim (B)$.
It is clear that the dimensions of irreducible representations
of $A \oplus_C B$ can be no larger than for $A$ and $B$.
There is a similar (nontopological) partition of
$\Prim_n (A \oplus_C B)$
as the disjoint union of homeomorphic copies of
$\Prim_n (\Ker (\rh)) \subset \Prim_n (A)$ and $\Prim_n (B)$.
Everything is metrizable.
The induced \ct\  field $E$ on $\Prim_n (A \oplus_C B)$ restricts
to the induced \ct\  fields on $\Prim_n (\Ker (\rh))$ and $\Prim_n (B)$,
which have finite type, so $E$ has finite type by Lemma~\ref{N12}.
Also,
\[
\dim (\Prim_n (A \oplus_C B))
   \leq \max (\dim (\Prim_n (A)), \, \dim (\Prim_n (B)) )
\]
by the same reasoning as in the last paragraph of the proofs of
Theorems~\ref{N13} and~\ref{N14}.
\end{pff}

\medskip

We will see in Example~\ref{N19}
below that subalgebras of separable \rsha s need not be \rsha s.
However, unitized ideals are.

\begin{cor}\label{N17.5}
Let $A$ be a separable \rsha, and let $I \subset A$ be an ideal.
Then the unitization $I^+$ is a \rsha,
with \tdim\  and \mms\  no larger than for $A$.
\end{cor}

\begin{pff}
We can write $I^+ = \C \oplus_{A/I} A$ using obvious maps,
and the result then follows from Proposition~\ref{N16}.
\end{pff}

\begin{prp}\label{N17}
(Compare with Theorem~11.8 of \cite{Pd}.)
Let $A$ and $B$ be separable \rsha s.
Then $A \otimes B$ is a \rsha, with \tdim\  at most the sum of the
\tdim s of $A$ and of $B$,
and with \mms\  at most the product of the \mms s of $A$ and of $B$.
\end{prp}

\begin{pff}
The proof is by induction on the length of a decomposition of $B$.
If $B = C (X, M_n)$, then the result is immediate from
Lemmas~\ref{A7} and~\ref{A8.5}.

Suppose therefore that the conclusion holds for some $B$, and consider
\[
A \otimes
 {\ts{ \left( B \oplus_{C ( X^{(0)}, \, M_n )}  C (X, M_n) \right) }}
 \cong (A \otimes B) \oplus_{A \otimes C ( X^{(0)}, \, M_n )}
                           (A \otimes  C (X, M_n)).
\]
Then $A \otimes B$ and $A \otimes  C (X, M_n)$ are \rsha s
satisfying the required bounds on the \tdim\  and \mms,  by the
induction hypothesis and by the initial case of the induction
respectively.
Corollary~\ref{N16} therefore implies that
$A \otimes \left( B \oplus_{C ( X^{(0)}, \, M_n )}  C (X, M_n) \right)$
is a \rsha\  satisfying the same bounds.
\end{pff}

\medskip

We now give an example of a separable unital \ca\   whose
irreducible representations all have dimension at most $2$,
but is which not a \rsha.

\begin{exa}\label{N18}
{}From Example~4.6 of \cite{Ph3} and the discussion following it,
we obtain a complex line bundle $L$, not of finite type,
over the second countable
locally compact space $X = \coprod_{n = 1}^{\infty} \C P^n$.
Let $E = L \oplus (X \times \C)$, and let $A$ be the unitization
$\Gm_0 (L (E))^+$.
Clearly $A$ has a single one dimensional representation,
and its other irreducible representations all have dimension $2$.
Moreover, $\Prim_2 (A) = X$ and the induced \ct\  field is $L (E)$.
As a vector bundle, $L (E) \cong E \otimes E^*$, so that it contains
$L \otimes (X \times \C)^* \cong L$ as a subbundle.
Since $L$ does not have finite type, neither does $L (E)$.
(Use Lemma~\ref{N5}~(2).)
So $L (E)$ does not have finite type as an $M_2$-bundle, by
Proposition~\ref{N8}.
Therefore $A$ does not have a \rshd,  by Theorem~\ref{N13}.
\end{exa}

It is perhaps interesting to point out that this example is a
direct limit of \rsha s, namely
$\Gm_0 (L (E))^+ \cong \dirlim \Gm (L (E |_{X_n}))^+$,
with $X_n = \coprod_{k = 1}^n \C P^k \subset X$.

\begin{exa}\label{N19}
Subalgebras of separable \rsha s, even of separable homogeneous \ca s,
need not be \rsha s.
We show this by constructing a unital embedding of the algebra in
Example~\ref{N18} above in a separable homogeneous \ca.
Let $Y$ be the Cantor set.
Choose a \ct\  surjective function $f$ from $Y$ to the one point
compactification $X^+$ of $X$.
(See, for example, Problem~O part~(e) in Chapter~5 of \cite{Kl}.)
Set $Y_n = f^{-1} (\C P^n )$, so that the sets $Y_n$ are disjoint
compact subsets of $Y$, and let $f_n = f |_{Y_n}$.
Then $\Gm (L (E |_{\C P^n}))$ can be canonically identified with the
set of sections $a$ of $L (f_n^* (E |_{\C P^n}))$ satisfying
$a (y_1) = a (y_2)$ whenever $f_n (y_1) = f_n (y_2)$.
(Recall that $f_n^* (E |_{\C P^n})_{y_1} = f_n^* (E |_{\C P^n})_{y_2}$
when $f_n (y_1) = f_n (y_2)$.)
This identification extends in an obvious way to $\Gm_0 (L (E))$.
However, $F = \left( f |_{f^{-1} (X)} \right)^* (E)$ is trivial,
because it is a bundle over the totally disconnected space $f^{-1} (X)$.
Therefore we obtain unital inclusions
\[
\Gm_0 (L (E))^+ \subset \Gm_0 (L (F))^+ \subset C (Y, M_2).
\]
\end{exa}

\section{Cancellation in \rsha s}

There are three basic theorems related to cancellation and subbundles of
vector bundles on finite dimensional compact spaces.
Roughly, they are that a vector bundle of sufficiently large rank
contains a trivial summand, that two stably isomorphic vector bundles
of sufficiently large rank are actually isomorphic, and that if the
difference between the ranks of two vector bundles
is sufficiently large, then the smaller one is a direct summand in
the bigger one.
In each case, ``sufficiently large'' means at least about half the
dimension of the space.
See Section~9.1 of \cite{Hs}, and see
Theorem~2.5 of \cite{Gd}, Lemma~3.4 of \cite{MP}, and
Lemma~1.5 of \cite{Ph1}
for restatements in terms of \pj s in the \ca s $C (X, M_n)$.
(Note that \cite{Hs} contains a slightly weaker version of the third
result, and that in the context of \cite{Hs} the third result
as stated above essentially contains the first.)
In this section, we generalize the second and third results to
\rsha s of finite \tdim,
and also prove an analogous result for unitaries.
By applying the characterization results of Section~2, we then
obtain results for type~1 \ca s whose irreducible representations have
bounded dimension.
In \cite{PhX}, these results will be applied to cancellation and
related problems in direct limits of \rsha s.

The first result fails, because there might not be any trivial
bundles of small rank.
Indeed, the \rsha s in the direct system implicit in Example~4.8
of \cite{PhX} have arbitrarily large minimum matrix size
but no nontrivial \pj s.
To prove the other two results, we first prove relative versions for
vector bundles over finite complexes.
Next, we switch to the \ca\  context (that is, $C (X, M_n)$)
and simultaneously generalize to arbitrary compact spaces.
The proofs for \rsha s can then be done by induction.
The results for unitaries can be gotten from the same preliminary
results as those for \pj s.

In this section, we use the
notation $p \sim q$ for Murray-von Neumann equivalence, and
$p \precsim q$ to mean that $p$ is Murray-von Neumann equivalent to a
sub\pj\  of $q$.
We let $U (A)$ denote the unitary group of a unital \ca\ $A$, and
we let  $U_0 (A)$ be the identity component of $U (A)$.
We write (as usual) $\GL_n (\C)$ for the invertible group $\inv (M_n)$.

Parts~(1) and~(2) of the first proposition
are relative versions of Theorems~9.1.2 and
9.1.5 of \cite{Hs}, except that part~(2) is generalized
in the manner of Theorem~2.5 of \cite{Gd}, Lemma~3.4 of \cite{MP}, and
Lemma~1.5 of \cite{Ph1}.
In part~(2), the homotopy condition in the hypotheses is necessary.
(For $d$ even, take $X$ to be the closed unit ball in $\R^d$,
take $Y$ to be its boundary $S^{d - 1}$, take $F = 0$,
take $E_1 = E_2$ to be trivial and suffiently large,
and take $a_0$ to represent a nontrivial element of $K^1 (S^{d - 1})$.)
The homotopy condition in the conclusion is then necessary for the
proposition to be useful in induction arguments.

\begin{prp}\label{B1}
Let $Y$ be a connected \cpt\  space, and let $X$ be a
connected \cpt\  space obtained from $X$ by attaching finitely
many cells of dimension at most $d$.
(In particular, the pair $(X, Y)$ could be a relative
CW-complex of dimension at most $d$.)

(1) Let $E$ and $F$ be (complex) vector bundles over $X$
with ranks satisfying $\rank (E) - \rank (F) \geq \frac{1}{2} (d - 1)$.
Let $a_0 \colon F |_Y \to E |_Y$ be an isomorphism of
$F |_Y$ with a subbundle of $E |_Y$.
Then there exists an isomorphism $a$ of
$F$ with a subbundle of $E$ such that $a |_Y = a_0$.

(2) Let $E_1$, $E_2$, and $F$ be vector bundles over $X$
with $\rank (E_1) \geq \frac{1}{2} d$.
Let $b \colon E_1 \oplus F \to E_2 \oplus F$ be an isomorphism,
let $a_0 \colon E_1 |_Y \to E_2 |_Y$ be an isomorphism,
and let $t \mapsto c_t^{(0)}$ be a homotopy of isomorphisms
$c_t^{(0)} \colon
   \left( E_1 \oplus F \right) |_Y \to \left( E_2 \oplus F \right) |_Y$
such that $c_0^{(0)} = a_0 \oplus \id_F |_Y$ and
$c_1^{(0)} = b |_Y$.
Then there exists an isomorphism $a \colon E_1 \to E_2$
such that $a |_Y = a_0$, and a homotopy of isomorphisms
$c_t \colon E_1 \oplus F \to E_2 \oplus F$ such that
\[
c_0 = a \oplus \id_{F}, \,\,\,\,\,\, c_1 = b,
\andeqn c_t |_Y = c_t^{(0)}.
\]
\end{prp}

\begin{pff}
By an obvious induction, we may assume $X$ is obtained from
$Y$ by attaching a single cell of dimension $n \leq d$, that
is, $X = Y \cup_f D^n$, where $D^n$ is the closed unit ball in
$\R^n$, and where $f$ is a \ct\  function from the boundary
$S^{n - 1}$ of $D^n$ to $Y$.
Note that $f$ induces a \ct\  function $g \colon D^n \to X$ which
agrees with $f$ on $S^{n - 1}$ and is a homeomorphism from
$D^n \setminus S^{n - 1}$ to $X \setminus Y$.
Pulling everything back via $f$ or $g$ as appropriate,
we further reduce to the case $Y = S^{n - 1}$, $X = D^n$,
and the attaching map is the identity.
Since $D^n$ is contractible, we may assume all bundles are
trivial.

For part~(1), let $\rank (E) = r$, so that $E = D^n \times \C^r$.
First suppose $\rank (F) = 1$, so that $F = D^n \times \C$.
In this case, $a_0$ is determined by a \ct\  function (section)
$\xi_0 \colon S^{n - 1} \to \C^r \setminus \{0\}$.
Now $\C^r \setminus \{0\}$ is homotopy equivalent to $S^{2 r - 1}$,
and the dimension hypotheses imply that $n - 1 < 2 r - 1$,
so $\xi_0$ extends to a \ct\  function (section)
$\xi \colon D^n \to \C^r \setminus \{0\}$.
This section gives the required isomorphism of $F$ with a
subbundle of $E$.

We now do the general case of part~(1) by induction on
$s = \rank (F)$.
Suppose the result is known for $F$ of rank $s$, and the actual
rank of $F$ is $s + 1$.
Write $F = F_0 \oplus (D^n \times \C)$.
Use the induction hypothesis to extend
$a_0 |_{\left( F_0 |_{S^{n - 1}} \right)}$ to an
isomorphism $c$ of $F_0$ with a subbundle of $E$ such that
$c |_{S^{n - 1}} = a_0 |_{\left( F_0 |_{S^{n - 1}} \right)}$.
Next, use the rank one case to extend
$a_0 |_{S^{n - 1} \times \C}$ to an
isomorphism $d$ of $D^n \times \C$ with a subbundle of the
orthogonal complement in $E$ of $c (F_0)$, such that
$d |_{S^{n - 1}} = a_0 |_{S^{n - 1} \times \C}$.
Then set $a = c + d$.

Next we do~(2).
\Wolog\  $E_1 = E_2 = D^n \times \C^r$ and $F = D^n \times \C^s$
for suitable $r$ and $s$.
The given maps can then be thought of as
$a_0 \colon S^{n - 1} \to \GL_r (\C)$, $b \colon D^n \to \GL_{r + s} (\C)$,
and $c_t^{(0)} \colon S^{n - 1} \to  \GL_{r + s} (\C)$,
such that  $t \mapsto c_t^{(0)}$ is a homotopy
from $a_0 \oplus 1$ to $b |_{S^{n - 1}}$.
In particular, $a_0$ defines a class in $\pi_{n - 1} ( \GL_r (\C))$
whose image in $\pi_{n - 1} ( \GL_{r + s} (\C))$ is zero.
Using the fact that the unitary groups are deformation retracts
of the corresponding invertible groups, Theorem~8.4.1 of \cite{Hs}
implies that
$\pi_{n - 1} ( \GL_r (\C)) \to \pi_{n - 1} ( \GL_{r + s} (\C))$
is an isomorphism when $n - 1 \leq 2 r - 1$.
This inequality is satisfied in our case, so $a_0$ is null homotopic.
Therefore there is a \ct\  function $\widetilde{a} \colon D^n \to \GL_r (\C)$
which extends $a_0$.

Let $h \colon [0, 1] \times D^n \to D^{n + 1}$ be a homeomorphism.
Then, with
\[
S = \left(\{ 0, 1\} \times D^n \right) \cup
                     \left([0, 1] \times S^{n - 1} \right),
\]
we have $h (S) = S^n$.
Let $\widetilde{v} \colon S^n \to  \GL_{r + s} (\C)$
be the function corresponding to the
function on $S$ given by
$\widetilde{a} (x) \oplus 1$ on $\{ 0 \} \times D^n$,
by $b (x)$ on $\{ 1 \} \times D^n$, and
by $c_t^{(0)} (x)$ on $[0, 1] \times S^{n - 1}$.
Let $\et = [\widetilde{v}] \in \pi_{n} ( \GL_{r + s} (\C))$.
Using Theorem~8.4.1 of \cite{Hs} in the same manner as above,
we find that
$\pi_n ( \GL_r (\C)) \to \pi_n ( \GL_{r + s} (\C))$ is surjective.
So there is a function $w \colon S \to \GL_r (\C)$ such that
$[w] \in \pi_n ( \GL_r (\C))$ has image
$- \et \in \pi_n ( \GL_{r + s} (\C))$.
A simple deformation argument allows us to require in addition that
$w (x) = 1$ for
$x \in \left( \{ 1\} \times D^n \right) \cup
                    \left( [0, 1] \times S^{n - 1} \right)$.
We have
$[ w \oplus 1] + [\widetilde{v}] = [(w \oplus 1) \widetilde{v}]$
in $\pi_n ( \GL_{r + s} (\C))$ (see Corollary 1.6.10 of \cite{Sp}),
whence $(w \oplus 1) \widetilde{v}$ is null homotopic.
Define $a = w \widetilde{a}$.
Since $w = 1$ on $\{ 0\} \times  S^{n - 1}$, we also have
$a |_{S^{n - 1}} = a_0$.
Further define $v \colon S \to \GL_{r + s} (\C)$ by substituting
$a$ for $\widetilde{a}$ in the definition of $\widetilde{v}$.
Then $v = (w \oplus 1) \widetilde{v}$ is null homotopic,
and hence extends continuously to a function
$(t, x) \mapsto c_t (x)$ from $[0, 1] \times D^n$ to
$\GL_{r + s} (\C)$.
Note that
\[
c_t |_{S^{n - 1}} = v |_{ \{ t\} \times S^{n - 1}}
 = \widetilde{v} |_{ \{ t\} \times S^{n - 1}} = c_t^{(0)} |_{S^{n - 1}}.
\]
So $t \mapsto c_t$ is the required homotopy.
\end{pff}

\begin{prp}\label{B3}
Let $X$ be a \cpt\  space with covering dimension $\dim (X) \leq d$,
and let $Y  \subset X$ be closed.

(1) Let $p, \, q \in C (X, M_n)$ be \pj s.
Suppose $\rank (p (x)) - \rank (q(x)) \geq \frac{1}{2} (d - 1)$
for all $x \in X$.
Let $s_0 \in C (Y, M_n)$ satisfy $s_0^* s_0 = q |_Y$ and
$s_0 s_0^* \leq p |_Y$.
Then there is $s \in C (X, M_n)$ such that
\[
s^* s = q, \,\,\,\,\,\, s s^* \leq p, \andeqn s |_Y = s_0.
\]

(2) Let $p_1, \, p_2, \, q_1, \, q_2 \in C (X, M_n)$ be \pj s,
and let $u, \, s \in C (X, M_n)$. Assume that
$\rank (p_j (x)) \geq \frac{1}{2} d$ for all $x$, that $p_j \perp q_j$, and that
\[
s^* s = q_1, \,\,\,\,\,\, s s^* = q_2, \,\,\,\,\,\,
u^* u = p_1 + q_1, \andeqn u u^* = p_2 + q_2.
\]
Further let $v_0 \in C (Y, M_n)$ satisfy
$v_0^* v_0 = p_1 |_Y$ and $v_0  v_0^* = p_2 |_Y$,
and let $t \mapsto w_t^{(0)}$ be a \ct\  path of partial isometries
in $C (Y, M_n)$ such that
\[
{\ts{ \left( \rsz{ w_t^{(0)} } \right)^* }}  w_t^{(0)}
      = \left( p_1 + q_1 \right) |_Y
                    \andeqn
w_t^{(0)} {\ts{ \left( \rsz{ w_t^{(0)} } \right)^* }}
      = \left( p_2 + q_2 \right) |_Y,
\]
and
\[
w_0^{(0)} = u |_Y \andeqn   w_1^{(0)} = v_0 + s |_Y.
\]
Then there is $v \in C (X, M_n)$ such that
\[
v^* v = p_1, \,\,\,\,\,\, v v^* = p_2, \andeqn v |_Y = v_0,
\]
and a \ct\  path $t \mapsto w_t$ of partial isometries
in $C (X, M_n)$ such that
\[
w_t^*  w_t = p_1 + q_1, \,\,\,\,\,\, w_t w_t^* = p_2 + q_2,
          \,\,\,\,\,\, w_0 = u, \,\,\,\,\,\, w_1 = v + s,
             \andeqn  w_t |_Y =  w_t^{(0)}.
\]
\end{prp}

\begin{pff}
We reduce this to Proposition~\ref{B1} in a number of steps.

The first reduction to to the case that all the \pj s involved have
constant ranks.
In~(1), this means
that $\rank (p (x))$ and $\rank (q (x))$ are independent of $x;$
similarly for~(2).
This is accomplished by writing $X$ as a finite disjoint union of
closed and open subsets on which all the ranks are constant.

Second, we reduce part~(2) to the situation $q_1 = q_2 = s = q$.
To do this, first replace (without renaming)
$M_n$ by $M_{2n}$, and replace each of
$p_j, \, q_j, \, u, \, s, \, v_0$, and $t \mapsto w_t^{(0)}$ by
its block diagonal direct sum with the $n \times n$ zero matrix,
namely $p_j \oplus 0$, $q_j \oplus 0$, etc.
Validity of the conclusion in the new situation is equivalent to its
validity in the old.
Now set
\[
z = \left( \begin{array}{cc}  s        & 1 - s s^*  \\
                            1 - s^* s  &    s^*     \end{array} \right),
\]
which is a unitary satisfying $z (q_1 \oplus 0) z^* = q_2 \oplus 0$.
Define
\[
\widetilde{q} = q_2 \oplus 0 = z (q_1 \oplus 0) z^* = (s \oplus 0) z^*,
           \,\,\,\,\,\,
\widetilde{p}_1 = z (p_1 \oplus 0) z^*,   \andeqn
\widetilde{p}_2 = p_2 \oplus 0,
\]
and
\[
\widetilde{u} = (u \oplus 0) z^*,  \,\,\,\,\,\,
\widetilde{v}_0 = (v_0 \oplus 0) z^*,  \andeqn
\widetilde{w}_t^{(0)}
   = {\ts{ \left( \rsz{ w_t^{(0)} \oplus 0 } \right) }} z^*.
\]
Assuming the proposition holds for this case (with
$\widetilde{q}_1 = \widetilde{q}_2 = \widetilde{s} = \widetilde{q}$),
let $\widetilde{v}$
and $\widetilde{w}_t$ be the resulting unitary and homotopy.
Then the conclusion~(2) as stated holds with
\[
v =
\left( \begin{array}{cc} 1 & 0 \\
                         0 & 0  \end{array} \right)
 \widetilde{v} z
\left( \begin{array}{cc} 1 & 0 \\
                         0 & 0  \end{array} \right)      
   \andeqn
w_t =
\left( \begin{array}{cc} 1 & 0 \\
                         0 & 0  \end{array} \right)
 \widetilde{w}_t z
\left( \begin{array}{cc} 1 & 0 \\
                         0 & 0  \end{array} \right),
\]
regarded as elements of $C (X, M_n)$.

In preparation for our next reduction, we prove the following
approximation result.
There is an absolute constant $\ep > 0$ such that, if the hypotheses
of~(1) are satisfied, and there are a \cpt\  space $\widetilde{X}$,
a closed subset $\widetilde{Y} \subset \widetilde{X}$,
a \ct\  function
$h \colon X \to \widetilde{X}$ with $h (Y) = \widetilde{Y}$, and
$\widetilde{p}, \, \widetilde{q}
                     \in C \left( \rsz{ \widetilde{X}, M_n } \right)$
and $\widetilde{s}_0 \in C \left( \rsz{ \widetilde{Y}, M_n } \right)$
for which the hypotheses and conclusion hold, moreover with
\[
\| \widetilde{p} \circ h - p \| < \ep,  \,\,\,\,\,\,
\| \widetilde{q} \circ h - q \| < \ep,  \andeqn
\| \widetilde{s}_0 \circ h - s_0 \| < \ep,
\]
then the conclusion of~(1) holds for $p, \, q$ and $s_0$.
The analogous statement is true for part~(2), in the case
$q_1 = q_2 = s = q$.

To see this in part~(1), let
$\widetilde{s} \in C \left( \rsz{ \widetilde{X}, M_n } \right)$
be the partial isometry obtained from the assumption that~(1) holds
for $\widetilde{X}$ and $\widetilde{Y}$.
Let $a \in C (X, M_n)$ be an arbitrary extension of $s_0$.
Choose a \nbhd\ $U$ of $Y$ such that
$\| \widetilde{s} \circ h (x) - a (x) \| < 2 \ep$ for $x \in U$.
(This norm is less than $\ep$ on $Y$.)
Choose a \ct\  function $f \colon X \to [0, 1]$ which is equal to $1$
on $Y$ and vanishes outside $U$.
Define
\[
b (x) = p (x)
 \left[ \rule{0em}{2ex} f (x) a (x)
                   + (1 - f (x)) \widetilde{s} \circ h (x) \right]
                                                       q (x)
                   \andeqn
s (x) = b (x) \left[ b (x)^* b (x) \right]^{-1/2}
\]
(functional calculus in $q (x) M_n q (x)$).
It is easily seen that if $\ep$ is small enough, then $s$ is the
required partial isometry.

The proof for~(2) is similar.
The formulas are as follows.
Let $a \in C (X, M_n)$ be an arbitrary extension of $v_0$.
Let $f \colon X \to [0, 1]$ be equal to $1$ on $Y$ and vanish outside
a suitable \nbhd\  of $Y$, and set
\[
b (x) = p_2 (x)
 \left[ \rule{0em}{2ex} f (x) a (x)
                   + (1 - f (x)) \widetilde{v} \circ h (x) \right]
                                                       p_1 (x)
                   \andeqn
v (x) = b (x) \left[ b (x)^* b (x) \right]^{-1/2}
\]
(functional calculus in $p_1 (x) M_n p_1 (x)$).
Then let $(t, x) \mapsto c (t, x)$ be an element of
$C ([0, 1] \times X, \, M_n)$ which extends the function
\[
(t, x) \mapsto \left\{ \begin{array}{lll}
         w_t^{(0)} (x)   & &  t \in [0, 1], \,\, x \in Y  \\
         v (x) + q (x)   & &  t = 1, \,\, x \in X  \\
         1               & &  t = 0, \,\, x \in X.
                                           \end{array} \right.
\]
Let $g \colon [0, 1] \times X \to [0, 1]$ be equal to $1$ on
$([0, 1] \times Y) \cup (\{0, 1\} \times X)$ and vanish outside
a suitable \nbhd\  of this set, and set
\[
d (t, x) = [p_2 (x) + q (x)]
   \left[ \rule{0em}{2ex} g (t, x) c (t, x)
                  + (1 - g (t, x)) \widetilde{w}_t \circ h (x) \right]
                                                  [p_1 (x) + q (x)]
\]
and
\[
w_t (x) = d (t, x) \left[ d (t, x)^* d (t, x) \right]^{-1/2}
\]
(functional calculus in $[p_1 (x) + q (x)] M_n [p_1 (x) + q (x)]$).

It is clear from this approximation result
that if $X = \invlim X_{\af}$, and the
proposition holds for all pairs $(X_{\af}, Y)$ with $Y \subset X_{\af}$
closed, then it holds for all pairs $(X, Y)$ with $Y \subset X$ closed.
(The function $h$ will be the map to a suitable $X_{\af}$.)
Our third reduction is from $X$ a compact space of dimension at most
$d$ to $X$ a compact {\emph{metric}} space of dimension at most $d$.
This follows because $\dim (X) \leq d$ implies $X \cong \invlim X_{\af}$
with $X_{\af}$ compact metric and $\dim ( X_{\af}) \leq d$ (Theorem
3.3.7 of \cite{En}).
For the fourth reduction, apply Theorem~1.13.5 of \cite{En}
(every compact space of dimension at most $d$ is the inverse limit
of finite simplicial complexes dimension at most $d$) to reduce to
the case in which $X$ is a finite simplicial complex (but still $Y$
is an arbitrary closed subset of $X$).

Our fifth reduction is to the situation in which $Y$ is a subcomplex.
An argument very similar to the approximation argument above allows
us to assume that $s_0$ in part~(1), and $v_0$ and $w_t^{(0)}$ in
part~(2), are actually defined on an open set $U \supset Y$.
Let $\ep = \dist (Y, X \setminus U) > 0$.
By repeated barycentric subdivision, we may assume that all the
simplexes in $X$ have diameter less than $\frac{1}{2} \ep$.
(See the proofs of Corollary~3.3.13 and Theorem~3.3.14 of \cite{Sp}.)
Let $\widetilde{Y}$ be the closed subset of $X$ given as the union of
all faces (of any dimension) of all simplexes in $X$ which intersect
$Y$.
The conclusion is now proved by applying the proposition for the case
of a subcomplex to $s_0 |_{\widetilde{Y}}$ in part~(1), and to
$v_0 |_{\widetilde{Y}}$ and $w_t^{(0)} |_{\widetilde{Y}}$ in part~(2).

It remains only to prove the result for the case that $Y$ is a
subcomplex of $X$.
For~(1), apply Proposition~\ref{B1}~(1)
with $E$ and $F$ the vector bundles
determined by $E_x = p (x) \C^n$ and $F_x = q (x) \C^n$, and with
$a_0 = s_0$.
(This is Swan's Theorem; see Theorem~2 of \cite{Sw}.)
The resulting isomorphism $a$ of $F$ with a subbundle of $E$
is an element of $C (X, M_n)$ such that $a^* a$ is invertible in
$q C (X, M_n) q$, so we can set $s = a (a^* a)^{ - 1/2}$.
The same application of Swan's Theorem and unitarization reduces
part~(2) (as reduced in the first reduction of this proof) to part~(2)
of Proposition~\ref{B1}.
\end{pff}

\medskip

Now we can prove the analogs of Theorem~2.5~(b) and~(c) of \cite{Gd}.
(The analog of Theorem~2.5~(a) of \cite{Gd} does not make sense in this
context, because there might not be any rank one \pj s in a \rsha.)

We follow Definition~5.1.1 of \cite{Bl2}, and write
$\Mi (A)$ for the {\emph{algebraic}} direct limit
$\dirlim M_n (A)$ under the maps $a \mapsto a \oplus 0$.

In part~(2) of the following proposition, we would really like to
allow $\rank (\ev_x (e))$ to be either zero or
greater than $\frac{1}{2} d$.
We have not been able to decide whether this is possible;
there might be a topological obstruction.

\begin{prp}\label{B4}
Let $A$ be a \rsha\  with total space $X$, and with
topological dimension function $d \colon X \to \N \cup \{0\}$
(as in Definition~\ref{A2}).
Let $p, \, e, \, f, \, q \in \Mi (A)$ be \pj s.

(1) If for every $x \in X$, either
\[
q (x) = 0 \,\,\,\,\,\, {\text{or}} \,\,\,\,\,\,
\rank (\ev_x (p)) - \rank (\ev_x (q))
      \geq {\ts{ \frac{1}{2} }} [d (x) - 1],
\]
then $q \precsim p$.

(2) If $e \oplus q \sim f \oplus q$ and
$\rank (\ev_x (e)) \geq \frac{1}{2} d (x)$ for every $x \in X$,
then $e \sim f$.
\end{prp}

\begin{pff}
We may assume everything is in $M_n (A)$ for some $n$.
Since $M_n (A)$ is also a \rsha\  with the same base spaces
(see Lemma~\ref{A7}) and hence the same topological dimension function,
we need only consider \pj s in $A$.

For~(1), we first carry out a reduction.
Adopt the notation of Definition~\ref{A2}.
Let $q_k$ be the image of $q$ in $C (X_k, \, M_{n (k)} )$.
Define
\[
Y_k = \{ x \in X_k \colon q_k (x) = 0 \} \andeqn Z_k = Y_k \setminus X_k.
\]
Applying Lemma~\ref{A5.1}, we construct a new decomposition for $A$
using the spaces $Y_k$ and $Z_k$ instead of $X_k$.
It has length at most $2 l + 1$, the same total space, and the new
topological dimension function is dominated by the old one,
since $\dim (X_k) = \max ( \dim (Y_k), \, \dim (Z_k) )$.
Using the old notation for this new decomposition, it has the
property that, for every $k$, either $q_k = 0$ or
$\rank (q_k (x)) \geq \frac{1}{2} [d (x) - 1]$ for all $x \in X_k$.

We now prove the result by induction on the length $l$ of a decomposition
with this property.
If the length is zero, then $A = C(X, M_n)$ for some $X$ and $n$.
In this case, if $q = 0$ then the conclusion is trivial, while if
$\rank (q (x)) \geq \frac{1}{2} [\dim (X) - 1]$ for all $x \in X$
then the conclusion follows from the case
$Y = \varnothing$ of Proposition~\ref{B3}~(1).

Now suppose the result is known for length $l$.
Let $A = B \oplus_{C ( X^{(0)}, \, M_n )} C (X, M_n)$, with
$\ph \colon B \to C \left( \rsz{ X^{(0)} }, \, M_n \right)$ unital and
$\rh \colon C (X, M_n) \to C \left( \rsz{ X^{(0)} }, \, M_n \right)$ the restriction map,
where $B$ is a \rsha\  of length $l$.
Write $p = (p_1, p_2)$ and $q = (q_1, q_2)$ with $p_1, \, q_1 \in B$
and $p_2, \, q_2 \in C (X, M_n)$.
By the induction assumption, there is $s_1 \in B$ such that
$s_1^* s_1 = q_1$ and $s_1 s_1^* \leq p_1$.

If $q_2 = 0$, then $\ph (s_1) = 0$.
Therefore $s = (s_1, 0)$ satisfies $s^* s = q$ and $s s^* \leq p$,
so $q \precsim p$ as desired.
Otherwise, $\rank (q_2 (x)) \geq \frac{1}{2} [\dim (X) - 1]$
for all $x \in X$.
Apply Proposition~\ref{B3}~(1) to $p_2, \, q_2, \, X^{(0)} \subset X$,
and with the partial isometry $\ph (s_1) \in C \left( \rsz{ X^{(0)} }, \, M_n \right)$.
Let $s_2$ be the resulting partial isometry.
Then $s = (s_1, s_2)$ is a partial isometry implementing the
relation $q \precsim p$.
This completes the induction, and proves~(1).

For~(2), we also use induction on the length.
We actually prove a stronger result (needed for the induction to work):
given a partial isometry $u$ with $u^* u = e + q$ and $u u^* = f + q$,
there is a partial isometry $s$ with $s^* s = e$, $s s^* = f$, and
$s + q$ homotopic to $u$ among partial isometries from $e + q$ to
$f + q$.

In the initial step (length $0$),
the conclusion follows from the case
$Y = \varnothing$ of Proposition~\ref{B3}~(2).
In the induction step, use the same notation as in the induction step
for~(1), and in addition let $e = (e_1, e_2)$, $f = (f_1, f_2)$,
and $u = (u_1, u_2)$.
The induction assumption provides a partial isometry $s_1 \in B$
such that $s_1^* s_1 = e_1$, $s_1 s_1^* = f$, and
there is a homotopy $t \mapsto w_1^{(t)}$ from
$s_1 + q_1$ to $u_1$ in the set of partial isometries from $e_1 + q_1$
to $f_1 + q_1$.
Use Proposition~\ref{B3}~(2) on $X^{(0)} \subset X$ to find
$s_2 \in C (X, M_n)$ and $t \mapsto w_2^{(t)} \in C (X, M_n)$,
extending $\ph (s_1)$ and
$t \mapsto \ph {\ts{ \left( \rsz{ w_2^{(t)} } \right) }}$ respectively,
such that $s_2^* s_2 = e_2$, $s_2 s_2^* = f$, and
$t \mapsto w_2^{(t)}$ is a homotopy from
$s_2 + q_2$ to $u_2$ in the set of partial isometries from $e_2 + q_2$
to $f_2 + q_2$.
Define
\[
w^{(t)} = {\ts{ \left( \rsz{ w_1^{(t)}, \, w_2^{(t)} } \right) }}
\andeqn s = (s_1, s_2).
\]
Then $s$ is a partial isometry from $e$ to $f$, and
$t \mapsto w^{(t)}$ is a homotopy from
$s + q$ to $u$ in the set of partial isometries from $e + q$
to $f + q$.
This completes the induction, and the proof of~(2).
\end{pff}

\medskip

The following proposition is an analog of Proposition~\ref{B4}
for unitaries instead of \pj s.

\begin{prp}\label{B5}
Let $A$ be a \rsha\  with total space $X$, and with
topological dimension function $d \colon X \to \N \cup \{0\}$
(as in Definition~\ref{A2}).
Let $p, \, q \in \Mi (A)$ be \pj s with $p \leq q$.

(1) If $\rank (\ev_x (p)) \geq \frac{1}{2} d (x)$ for all $x \in X$,
and $u \in U (q \Mi (A) q)$, then there exists $v \in U (p \Mi (A) p)$
such that $v + (q - p)$ is homotopic to $u$ in $U (q \Mi (A) q)$.

(2) If $\rank (\ev_x (p)) \geq \frac{1}{2} [d (x) + 1]$ for all $x \in X$,
and $v^{(0)}, \, v^{(1)} \in U (p \Mi (A) p)$ are unitaries such that
$v^{(0)} + (q - p)$ is homotopic to $v^{(1)} + (q - p)$ in
$U (q \Mi (A) q)$,
then $v^{(0)}$ is homotopic to $v^{(1)}$ in $U (p \Mi (A) p)$.
\end{prp}

\begin{pff}
As in the proof of Proposition~\ref{B4}, we may assume $q \in A$.

We prove the result by induction on the length $l$ of a decomposition.
We describe only the induction steps, since the initial steps differ
only in that the subset $Y$ in Proposition~\ref{B3} is taken to be empty
(in~(1)) or just $\{0, 1\} \times X$ (in~(2)).

So suppose part~(1) is known for \rsha s of length $l$.
Let $A = B \oplus_{C ( X^{(0)}, \, M_n )} C (X, M_n)$, with
$\ph \colon B \to C \left( \rsz{ X^{(0)} }, \, M_n \right)$ unital and
$\rh \colon C (X, M_n) \to C \left( \rsz{ X^{(0)} }, \, M_n \right)$ the restriction map,
where $B$ is a \rsha\  of length $l$.
Write $p = (p_1, p_2)$, $q = (q_1, q_2)$, and $u = (u_1, u_2)$,
with $p_1, \, q_1, \, u_1 \in B$
and $p_2, \, q_2, \, u_2 \in C (X, M_n)$.
By the induction assumption, there are $v_1 \in U (p_1 B p_1)$
and $t \mapsto w_1^{(t)} \in U (q_1 B q_1)$ such that
$w_1^{(0)} = u_1$ and $w_1^{(1)} = v_1 + (q_1 - p_1)$.
Apply Proposition~\ref{B3}~(2) to $X^{(0)} \subset X$, with the
$p_i$ there both being $p_2$, the $q_i$ both being $q_2 - p_2$,
$v_0$ being $\ph (v_1)$, $s$ being $q_2 - p_1$, and the homotopy
being $t \mapsto \ph {\ts{ \left( \rsz{ w_1^{(t)} } \right) }}$.
Let $v_2$ be the resulting element of $U (p_2  C (X, M_n) p_2)$
and let $t \mapsto w_2^{(t)} \in U (q_2  C (X, M_n) q_2)$ be the
resulting homotopy.
The elements $v = (v_1, v_2)$ and
$w^{(t)} = {\ts{ \left( \rsz{ w_1^{(t)}, w_2^{(t)} } \right) }}$
prove the induction step.

For part~(2), we use analogous notation, and we assume for the
induction step that we are given homotopies $t \mapsto c^{(t)}$
from $v^{(0)} + (q - p)$ to $v^{(1)} + (q - p)$ in
$U (q A q)$ and $t \mapsto v_1^{(t)}$ from $v_1^{(0)}$ to $v_1^{(1)}$
in $U (p_1 B p_1)$.
We regard the homotopies as elements of the \rsha\  $C ([0, 1], \, A)$
with its obvious corresponding decomposition (see Lemma~\ref{A8.5}),
and we require as part of the induction hypothesis
that there be a homotopy
$(t_1, t_2) \mapsto w_1^{(t_1, t_2)} \in U (q_1 B q_1)$
of elements of $C ([0, 1], \, B)$ such that
\[
w_1^{(0, t_2)} = v_1^{(t_2)} + (q_1 - p_1) \andeqn
w_1^{(1, t_2)} = c_1^{(t_2)}
\]
for all $t_2$, and
\[
w_1^{(t_1, 0)} = v_1^{(0)} + (q_1 - p_1) \andeqn
w_1^{(t_1, 1)} = v_1^{(1)} + (q_1 - p_1)
\]
for all $t_1$.
We apply Proposition~\ref{B3}~(2) to
\[
{\ts{ \left( \rsz{  [0, 1] \times X^{(0)} } \right) }}
                        \cup \left( \{0, 1\} \times X \right)
       \subset [0, 1] \times X,
\]
with the $p_i$ there both being $(t, x) \mapsto p_2 (x)$, with
the $q_i$ both being $(t, x) \mapsto q_2 (x) - p_2 (x)$, with
$v_0$ there being
$(t, x) \mapsto \ph {\ts{ \left( \rsz{ v_1^{(t)} } \right) }} (x)$ on
$[0, 1] \times X^{(0)}$ and $v_2^{(i)} (x)$ on $\{ i \} \times X$, with
$u$ there being $(t, x) \mapsto c_2^{(t)} (x)$, with
$s$ there being $(t, x) \mapsto q_2 (x) - p_2 (x)$,
and with $(t, y) \mapsto w_t^{(0)} (y)$ being
\[
(t_1, t_2, x) \mapsto \left\{
  \begin{array}{ll}
        \ph {\ts{ \left( \rsz{ w_1^{(t_1, t_2)} } \right) }} (x)
                                 & x \in X^{(0)} \\
           v_2^{(t_2)} (x) + (q_2 (x) - p_2 (x)) & t_2 \in \{0, 1\}.
                                                   \end{array} \right.
\]
(Here $t_1 = t$ and $(t_2, x) = y$.)
Note that $\dim ([0, 1] \times X) \leq \dim (X) + 1$ by
Proposition~3.2.6 of \cite{Pr}.
We then obtain a unitary $(t, x) \mapsto v_2^{(t)} (x)$ in
$q_2 [ C ([0, 1] \times X, \, M_n) ] q_2$ extending
$(t, x) \mapsto \ph {\ts{ \left( \rsz{ v_1^{(t)} } \right) }} (x)$
and agreeing at $t = 0$ and
$t = 1$ with $v_2^{(0)}$ and $v_2^{(1)}$ as already defined,
and also a homotopy
$(t_1, t_2) \mapsto w_2^{(t_1, t_2)} \in U (q_2 C ( X, M_n) q_2)$
from $t_2 \mapsto v_2^{(t_2)} + (q_2 - p_2)$ to
$t_2 \mapsto c_2^{(t_2)}$ extending the given homotopy,
in particular satisfying
$w_2^{(t_1, t_2)} (x) = v_2^{(t_2)} (x) + (q_2 (x) - p_2 (x))$
for $x \in X$, $t_1 \in [0, 1]$, and $t_2 \in \{0, 1\}$.
The induction step is completed by taking
\[
v^{(t)} = {\ts{ \left( \rsz{ v_1^{(t)}, v_2^{(t)} } \right) }} \andeqn
w^{(t_1, t_2)}
  = {\ts{ \left( \rsz{ w_1^{(t_1, t_2)}, w_2^{(t_1, t_2)} } \right) }}.
\]
One readily checks that the properties assumed in $B$ have been
extended to $A$.
\end{pff}

\medskip

Although the primary intended application of the results of this
section is to direct limits of \rsha s (see \cite{PhX}),
there are immediate interesting consequences for type~1 \ca s.
As before Theorem~\ref{N10.5},
we let $\Prim_n (A)$ denote the subspace of $\Prim (A)$ consisting
of the kernels of $n$-dimensional representations of $A$.

\begin{thm}\label{Y1}
Let $A$ be a separable unital \ca.
Suppose that there is $N \in \N$ such that all
irreducible representations of $A$ have dimension at most $N$.
Let $p, \, q \in \Mi (A)$ be \pj s.
Suppose that for every $n$ and every irreducible representation
$\pi$ of $A$ of dimension $n$, we have
\[
\rank (\pi (p)) - \rank (\pi (q)) \geq
  {\ts{ \frac{1}{2} }} [ \dim (\Prim_n (A)) - 1].
\]
Then $q \precsim p$.
\end{thm}

\begin{pff}
Lemma~\ref{N4} and Proposition~\ref{N8} allow us to apply
Proposition~\ref{N11}.
Using the standard notation (see Definition~\ref{A2}) for the
resulting \rshd, we then have:
\bit
\item[(1)]
$n (0) \leq n (1) \leq \cdots \leq n (l)$.
\item[(2)]
$\dim (X_k) \leq \dim \left( \Prim_{n (k)} (A) \right)$
for $0 \leq k \leq l$.
\item[(3)]
$X_k \setminus X_k^{(0)}$ is dense in $X_k$ for $1 \leq k \leq l$.
\eit
Let $X$ be the total space, and let $d \colon X \to \N \cup \{ 0 \}$ be
the \tdim\  function (as in Definition~\ref{A2}).
Take $X_0^{(0)} = \varnothing$.
By Lemma~\ref{A9}, the evaluations $\ev_x$,
for $x \in X_k \setminus X_k^{(0)}$, are irreducible.
Therefore, for  $x \in X_k \setminus X_k^{(0)}$, we have
\[
\rank (\ev_x (p)) - \rank (\ev_x (q)) \geq
  {\ts{ \frac{1}{2} }} [ \dim (\Prim_{n (k)} (A)) - 1] \geq
  {\ts{ \frac{1}{2} }} [ d (x) - 1].
\]
By continuity, this inequality is valid for all
$x \in {\overline{X_k \setminus X_k^{(0)} }} = X_k$.
Now Proposition~\ref{B4}~(1) implies that $q \precsim p$.
\end{pff}

\begin{thm}\label{Y2}
Let $A$ be a separable unital \ca.
Suppose that there is $N \in \N$ such that all
irreducible representations of $A$ have dimension at most $N$.
Let $e, \, f, \, q \in \Mi (A)$ be \pj s such that
$e \oplus q \sim f \oplus q$.
Suppose that for every $n$ and every irreducible representation
$\pi$ of $A$ of dimension $n$, we have
\[
\rank (\pi (e)) \geq {\ts{ \frac{1}{2} }} \dim (\Prim_n (A)).
\]
Then $e \sim f$.
\end{thm}

\begin{pff}
The proof is the same as that of Theorem~\ref{Y1},
except that in the last step we apply Proposition~\ref{B4}~(2).
\end{pff}

\begin{thm}\label{Y3}
Let $A$ be a separable unital \ca.
Suppose that there is $N \in \N$ such that all
irreducible representations of $A$ have dimension at most $N$,
and suppose that $\dim (\Prim_n (A)) \leq 2 n - 1$ for all $n$.
Then the natural map $U (A) / U_0 (A) \to K_1 (A)$ is an isomorphism.
\end{thm}

\begin{pff}
We first prove injectivity.
Thus, let $u, \, v \in A$ be unitary, and suppose that $[u] = [v]$
in $K_1 (A)$; we must show that $u$ can be connected to $v$ by a
path in $U (A)$.
By assumption, there is $n$ such that $u \oplus 1$
can be connected to $v \oplus 1$ in $U (M_n (A))$.
As in the proof of Theorem~\ref{Y1}, there is a \rshd\  for $A$ 
with total space $X$ and
topological dimension function $d \colon X \to \N \cup \{0\}$
such that $\rank (\ev_x (1_A)) \geq \frac{1}{2} [d (x) + 1]$
for every $x \in X$.
Apply Proposition~\ref{B5}~(2) with $p = 1_A$ and $q = 1_{M_n (A)}$.

Surjectivity is proved in the same way, using instead
Proposition~\ref{B5}~(1) at the end.
\end{pff}

\medskip

If $A = C (X, M_n)$ with $\dim (X) \leq 2 n - 1$, then one can
deduce the conclusion of Theorem~\ref{Y3}
from stable rank considerations.
Specifically, Proposition~1.7 of~\cite{Rf1} gives
$\tsr (C (X) ) \leq n$, and it then follows from
Theorem~2.10 of~\cite{Rf2} that
\[
U (C (X, M_n)) / U_0 (C (X, M_n)) \to K_1 (A)
\]
is an isomorphism.
However, it is not clear how to obtain the general case of
Theorem~\ref{Y3} from stable rank considerations.
As just one of several difficulties, we note that,
with $X$ and $n$ as above, Theorem~6.1 of~\cite{Rf1} gives
$\tsr (C (X, M_n)) = 2$ unless $\dim (X) \leq 1$.
Applying Theorem~4.3 of~\cite{Rf1}, we see that an algebra $A$
as in the hypotheses of Theorem~\ref{Y3} generally has stable
rank at least $2$, so that we can't apply Theorem~2.10 of~\cite{Rf2}.

\end{document}